\documentclass[[ejs,authoryear]{imsart}
\arxiv{math.ST: 1503.03212}

\startlocaldefs
\usepackage{amssymb,natbib}
\usepackage[cmex10]{amsmath}
\usepackage{amsmath,enumerate,amsbsy} 
\DeclareMathAlphabet{\mathcal}{OMS}{cmsy}{m}{n}
\newcommand{\bs}[1]{\boldsymbol{#1}}
 \usepackage{amsthm}
\newtheorem{theorem1}{Special Theorem}
\newtheorem{definition}[theorem1]{Definition}
\newtheorem{property}{Property}
\endlocaldefs

\begin{document}

\begin{frontmatter}
\title{Multivariate Generalized Gram-Charlier Series in Vector Notations}
\runtitle{Multivariate GGC Series}
\begin{aug}
\author{ \snm{  Dharmani Bhaveshkumar C.} 
\ead[label=e1]{dharmanibc@gmail.com }}
\address{Dhirubhai Ambani Institute of \\
       Information \& Communication Technology (DAIICT), \\
       Gandhinagar, Gujarat, INDIA - 382001 
      \printead{e1}
}
\runauthor{Dharmani Bhaveshkumar C.}
\end{aug}

 
\begin{abstract}
The article derives multivariate Generalized Gram-Charlier (GGC) series that expands an unknown joint probability density function (\textit{pdf}) of a random vector in terms of the differentiations of the joint \textit{pdf} of a reference random vector. Conventionally, the higher order differentiations of a multivariate \textit{pdf} in GGC series will require multi-element array or tensor representations. But, the current article derives the GGC series in vector notations. The required 
higher order differentiations of a multivariate \textit{pdf} in vector notations are achieved through application of a specific Kronecker product based differentiation operator. 
Overall, the article uses only elementary calculus of several variables; instead Tensor calculus; to achieve the extension of  an existing specific derivation for GGC series in univariate 
to multivariate. The derived multivariate GGC expression is more elementary as using vector notations compare to the  coordinatewise tensor notations and more comprehensive as apparently  more nearer to its counterpart for univariate. 
The same advantages are shared by the other expressions obtained in the article; such as the mutual relations between 
cumulants and moments of a random vector, integral form of a multivariate \textit{pdf}, integral form of the multivariate Hermite polynomials, the multivariate Gram-Charlier A (GCA) series and others. 
\end{abstract}
\end{frontmatter}

\section{Introduction} 
\label{introduction}
The Generalized Gram-Charlier (GGC) series expands an unknown \textit{pdf} 
as a linear combination of the increasing order differentiations of a reference \textit{pdf}, 
where the coefficients of expansion involve cumulant  differences between those of an unknown \textit{pdf} and a reference \textit{pdf}. 
The GGC expansions are used to approximate \textit{pdf}  and  functions of \textit{pdf} in Statistics, Machine Learning, 
Economics, Chemistry, Astronomy and other application areas.   
 There have been used Poisson's distribution \citep{aroian1937},  log-normal distribution \citep{GGC2pdf47}, binomial distribution \citep{bookReitz27}, gamma distribution \citep{bowers1966expansion} and others as reference \textit{pdf}s.
But, the Gram-Charlier (GC) expansion with Gaussian density as a reference \textit{pdf}
is the most popular and identified as the Gram-Charlier A (GCA) series. 
Specifically, the GCA series is used for test of Gaussianity or near Gaussian \textit{pdf}  approximations  \citep{GGC2pdf47,GGCchem94,girolaminegH96,GGCpdf01,GCDnd09}, for entropy measure and independence measure approximations  \citep{GGCHAmari96,ICAbook01}, for optima analysis through  derivatives of \textit{pdf} \citep{ICA03Boscolo}, for time-frequency analysis \citep{genGramC98} and others. The rearrangement of the terms in GCA series results into the Edgeworth series with better convergence property.   

\hspace{0.2 in} There exists multiple ways to derive univariate GCA series, as reported by \citet{HistoryGGS00,hald2002history}. It's generalization to univariate GGC  series is derived by \citet{Schleher77,genGramC98,genGramC07,genGramC11}.  


\hspace{0.2 in} The multivariate GGC or GCA series derivation requires multivariate representations of the Taylor series, the increasing  order differentiations of a reference \textit{pdf} and the cumulants. 
To signify the representation issue in required  multivariate extensions; it is worth quoting  
 \citet{ndHermite02} that says, ` $\ldots$ though the generalizations to higher dimensions may be  considered straightforward, the methodology and the mathematical notations get quite cumbersome  when dealing with the higher order derivatives of the characteristic function or other functions of a random  vector $\ldots$ '. 

\hspace{0.2 in}	Conventionally, the higher order differentiations of a multivariate \textit{pdf}; and therefore, the multivariate cumulants and multivariate Hermite polynomials; require multilinear representations. 
It is known and acknowledged historically that going from matrix like notations to tensor notations for multivariate cumulants and multivariate Hermite polynomials have made the representation more transparent and the calculations more simpler \citep{Tensorbook87}. 
Though the tensor notations have advantages over the matrix notations, they require componentwise separate representations and are more tedious compare to the vector notations. 
As a unified and comprehensive solution to this, there has been used an approach based only on elementary calculus of several variables by \citet{ndHermite02}.  The approach uses a specific Kronecker product based differential operator, identified as the `K-derivative operator',  to achieve vectorization of the Jacobian matrix of a multivariate vector  function. The successive applications of this K-derivative operator achieves vectorization of the higher order derivatives also. Using this approach,
there have been derived the multivariate Taylor series and the higher order cumulants  \citep{ndHermite02,SankhyaCumVec06}; as well,  the  multivariate Hermite polynomials \citep{ndHermite02} in vector notations resulting into more transparent  representations. In fact, it could be noticed that the same approach has been first used to derive vector Hermite polynomials by \citet{dnHermite96}; and later\footnote{Somehow, the citation for \citep{dnHermite96} is not found in article \citep{ndHermite02}.} it has been formalized and generalized by \citet{ndHermite02}. 

\hspace{0.2 in} There exists various approaches deriving multivariate GCA series with various representations. 
They include GCA series representation using multi-element matrix notations for moments and cumulants by \citet{sauer1979convenient}; using 
 recursive formula for Hermite polynomials by \citet{ndHPGCberkowitz70}; using tensor representation for cumulants and Hermite polynomials by \citet[Chapter 5]{Tensorbook87}; using tensor representation for cumulants and involving multivariate Bell polynomials by \citet{withers2014dual}; using vector moments and vector Hermite polynomials by \citet{dnHermite96} and others. 
  There exists various 
 derivations for multivariate Edgeworth series 
  \citep{ndEdge76,ndEdgeworth86,amari1983differential,Tensorbook87,unifieddensities98,withers2014dual}.  
  There also exists multivariate GGC series, derived by \citet[Chapter 5]{Tensorbook87}, in tensor notations. 
But, as per the author's knowledge, there exists neither the multivariate GGC series nor the multivariate Edgeworth series in vector notations. For the ease of the readers in following and comparing the various representations;
the existing representations of multivariate GGC series and multivariate GCA series are shown and compared in \ref{historyGGCreps}.

\hspace{0.2 in} Overall, to take advantages due to the recent advancement in representation, this article extends a  specific derivation for univariate GGC series by \citet{genGramC07} to multivariate; 
using only elementary calculus of several variables instead of Tensor calculus. 
As a by product, it  also derives mutual relations  between vector cumulants and vector moments of a random vector; integral form of the multivariate \textit{pdf}; integral form of the multivariate vector Hermite polynomials 
and the multivariate GCA series. All the derived multivariate expressions are more elementary as using vector notations  and more comprehensive as apparently  more nearer to their counterparts for univariate; compare to their coordinatewise tensor  notations.  The intermediate theoretical results, in the article, are verified using suitable known examples. 

\hspace{0.2 in} 
Towards the aim of the article, the next Section \ref{Kronprod} briefs some necessary background on the Kronecker product and a way to obtain vectorization of the higher order differentiations of a multivariate \textit{pdf}. It also obtains  the required multivariate Taylor series expansion using the derived notations.  
After the preliminary background, 
this article follows almost the same sequence for multivariate as that in \citep{genGramC07} for univariate. 
The Section \ref{chfn} uses the characteristic function and the generating functions to derive 
cumulants and moments of a random vector in vector notations with their mutual relationships. 
The Section \ref{pdfcum} obtains multivariate \textit{pdf} in terms of its vector cumulants. The expressions for derivatives of multivariate  Gaussian density and vector Hermite polynomials are derived in Section \ref{DerGauss}. The Section \ref{ndGCAseries} derives multivariate GCA series by representing an unknown \textit{pdf} in terms of the Gaussian \textit{pdf} as a reference. The Section \ref{ndGGC} derives GGC Expansion, representing an unknown \textit{pdf} in terms of the 
a known reference \textit{pdf}. The Section \ref{chGGC} derives an unknown characteristic function of a random vector in terms of a  reference characteristic function. The Section \ref{compactGGC} derives the same GGC expansion in a more compact way that summarizes  the approach of the whole derivation. Finally, Section \ref{conclusion} concludes the article. For the sake of clarity;  the calculation details, the proofs and the expressions for existing multivariate expansions are kept in appendix at the end of the article.     
\section{Vectorization of the higher order differentiations} 
\label{Kronprod}
The section briefs the Kronecker product and the way it can be applied to achieve vectorization of the higher order differentiations of a multivariate \textit{pdf}. Based on it, the multivariate Taylor series is obtained in vector notations. 
More details can be found 
on the Kronecker Product in \citep[Chapter 2]{MagnusBook99}, 
on achieving vectorization of the higher order differentiations in \citep{ndHermite02,SankhyaCumVec06} and on the commutation matrices in \citep[Chapter 3, Section 7]{MagnusBook99}.  
\begin{definition}[Kronecker Product Operator $(\otimes)$]
The Kronecker Product Operator $(\otimes)$ between matrices $\mathbf{A}$ with size $p \times q$ and $\mathbf{B}$ with size $m \times n$ is defined as:\\
\begin{align}
\mathbf{A}\otimes\mathbf{B} =  \left[
\begin{array}[]{c c c c}
a_{11}\mathbf{B} & a_{12}\mathbf{B} & \cdots & a_{1q}\mathbf{B} \\
a_{21}\mathbf{B} &a_{22}\mathbf{B} &\cdots &a_{2q}\mathbf{B} \\
\vdots & \vdots & \ddots&\vdots \\
a_{p1}\mathbf{B} &a_{p2}\mathbf{B} &\cdots &a_{pq}\mathbf{B}	
\end{array}
\right]
\end{align}
\end{definition}
The resultant matrix is of dimension $pm \times qn$. 
As a further example;  let $\mathbf{A}$ is with size $p \times 1$ and $\mathbf{B}$ is with size $m \times 1$, then $\mathbf{A} \otimes \mathbf{B}'$ is\footnote{The symbol ' stands for Transpose of a matrix} a matrix with size $p \times m $.  $\mathbf{A} \otimes \mathbf{A}$ is symbolically represented as $\mathbf{A}^{\otimes 2}$ and has size $p^2 \times 1$. In general, $\mathbf{A} \otimes \mathbf{A} \otimes \ldots \otimes \mathbf{A} \mbox{ (n times)}$ is symbolically represented as $\mathbf{A}^{\otimes n}$ and has size $p^n \times 1$.  
\begin{definition}[Jacobian Matrix] 
Let $\bs{\lambda} = (\lambda_1, \lambda_2, \ldots, \lambda_d)'$, $\bs{\lambda} \in \mathbb{R}^d$ and  $\mathbf{f}(\bs{\lambda}) = $ \\ 
$(f_1(\lambda), f_2(\lambda), \ldots, f_m(\lambda))' \in \mathbb{R}^{m}$ be a differentiable m-component vector function. Then, Jacobian matrix of $\mathbf{f}(\bs{\lambda})$  $(\mathbf{J}(\mathbf{f}))$ is an $m \times d$ matrix defined as under:
\begin{align}
\mathbf{J}(\mathbf{f}(\bs{\lambda})) = \frac{d \mathbf{f}}{d {\bs{\lambda}}} = \left[ \frac{\partial \mathbf{f}}{\partial \lambda_1}, \frac{\partial \mathbf{f}}{\partial \lambda_2},\ldots, \frac{\partial \mathbf{f}}{\partial \lambda_d} \right] = \left[
\begin{array}[]{c c c c}
\frac{\partial f_1}{\partial \lambda_1} & \frac{\partial f_1}{\partial \lambda_2}  & \cdots & \frac{\partial f_1}{\partial \lambda_d} \\
\frac{\partial f_2}{\partial \lambda_1} &  \ddots &  &\vdots  \\
\vdots  &  & \ddots  &\vdots \\
\frac{\partial f_m}{\partial \lambda_1} &  \frac{\partial f_m}{\partial \lambda_2} & \cdots & \frac{\partial f_m}{\partial \lambda_d} 
\end{array}
\right]
\end{align}
\end{definition}
Let the vector differential operator be defined as a column vector 
  $\mathbf{D}_{\bs{\lambda}} = \left( \frac{\partial}{\partial \lambda_1}, \frac{\partial}{\partial \lambda_2},\ldots, \frac{\partial}{\partial \lambda_d} \right)'$, then 
  the Jacobian matrix, in terms of the $\mathbf{D}_{\bs{\lambda}}$, can be re-written as: 
\begin{align}
\mathbf{J}(\mathbf{f}(\bs{\lambda})) = \mathbf{D}_{\bs{\lambda}}(\mathbf{f}) = \mathbf{f}(\bs{\lambda})\mathbf{D}'_{\bs{\lambda}}   = \left( f_1(\bs{\lambda}), f_2(\bs{\lambda}), \ldots, f_m(\bs{\lambda}) \right)' \left( \frac{\partial}{\partial \lambda_1}, \frac{\partial}{\partial \lambda_2},\ldots, \frac{\partial}{\partial \lambda_d} \right)  
\end{align}
This implies that to match the definition of differentiation from matrix calculus, the vector differential operator should be applied from the right to the left. This is same as the requirement to be satisfied on generalization of vector derivative to matrix derivative as discussed by \citet{Magnus10}. 
So, applying vector derivative operator from right to the left,  has been kept as a rule throughout the article. 
\begin{definition}[The K-derivative Operator] 
Let $\bs{\lambda} = (\lambda_1, \lambda_2, \ldots, \lambda_d)'$, $\bs{\lambda} \in \mathbb{R}^d$; the vector differential operator $\mathbf{D}_{\bs{\lambda}} = \left( \frac{\partial}{\partial \lambda_1}, \frac{\partial}{\partial \lambda_2},\ldots, \frac{\partial}{\partial \lambda_d} \right)'$ and  a differentiable m-component vector function 
$\mathbf{f}(\bs{\lambda}) = (f_1(\lambda), f_2(\lambda), \ldots, f_m(\lambda))' \in \mathbb{R}^{m}$.  Then, the  Kronecker product between $\mathbf{D}_{\bs{\lambda}}$ and $\mathbf{f}(\bs{\lambda})$ is given as under:
\begin{align}
\mathbf{D}_{\bs{\lambda}}^\otimes \mathbf{f}(\bs{\lambda}) &= 
\left[ 
\begin{array}[]{c }
f_1(\bs{\lambda}) \\ f_2(\bs{\lambda}) \\ \vdots \\ f_m(\bs{\lambda})
\end{array}
\right]
\otimes
\left[ 
\begin{array}[]{c }
\frac{\partial }{\partial \lambda_1} \\ \frac{\partial }{\partial \lambda_2} \\ \vdots \\ \frac{\partial }{\partial \lambda_d}
\end{array}
\right] =
Vec\left[
\begin{array}[]{c c c c}
\frac{\partial f_1}{\partial \lambda_1} & \frac{\partial f_1}{\partial \lambda_2}  & \cdots & \frac{\partial f_1}{\partial \lambda_d} \\
\frac{\partial f_2}{\partial \lambda_1} &  \ddots &  &\vdots  \\
\vdots  &  & \ddots  &\vdots \\
\frac{\partial f_m}{\partial \lambda_1} &  \frac{\partial f_m}{\partial \lambda_2} & \cdots & \frac{\partial f_m}{\partial \lambda_d} 
\end{array}
\right]' \\
\Rightarrow \mathbf{D}_{\bs{\lambda}}^\otimes \mathbf{f}(\bs{\lambda}) &= Vec\left(\frac{\partial \mathbf{f}}{\partial \bs{\lambda}'}\right)' = Vec \left( \frac{\partial}{\partial \bs{\lambda}}\mathbf{f}' \right)
\end{align}
where, the $Vec$ operator converts $m \times d$ matrix into an $ md \times 1$ column vector by stacking the columns one after an other. 
The operator $\mathbf{D}_{\bs{\lambda}}^{\otimes } $ is called Kronecker derivative operator  or simply, K-derivative operator. 
\end{definition} 
Thus, the Kronecker product with the vector differential operator, obtains vectorization of the transposed Jacobian of a vector function.  Corresponding to the definition, the $k^{th}$ order differentiation is given by:
\begin{align} 
 \mathbf{D}_{\bs{\lambda}}^{\otimes k} \mathbf{f} = \mathbf{D}_{\bs{\lambda}}^{\otimes}\left( \mathbf{D}_{\bs{\lambda}}^{\otimes k-1}\mathbf{f} \right) = [f_1(\lambda), f_2(\lambda), \ldots, f_m(\lambda)]' \otimes \left[ \frac{\partial }{\partial \lambda_1},  \frac{\partial }{\partial \lambda_2},  \cdots,  \frac{\partial }{\partial \lambda_d}  \right]^{'\otimes k} 
\end{align}
The $\mathbf{D}_{\bs{\lambda}}^{\otimes k}\mathbf{f}$ is a column vector of dimension $md^k \times 1$. 
Some important properties of the K-derivative operator, those are useful in the further derivations, are listed in Appendix \ref{KderProp}. 
\subsection{Application of the K-derivative operator to the multivariate Taylor series} 
Let $\mathbf{x} = \left(X_1, X_2,..., X_d\right)'$ be a d-dimensional column vector and $f(\mathbf{x})$ be the function of several variables differentiable in each variable. Using the defined K-derivative operator, the Taylor series for $f(\mathbf{x})$, expanding it at origin, is given as:
\begin{equation}
\label{ndTaylor}
f(\mathbf{x}) = \sum_{m=0}^{m=\infty}\frac{1}{m!}\mathbf{c}(m,d)'\mathbf{x}^{\otimes m}
\end{equation} 
where, $\mathbf{c}(m,d)$ is the vector of dimension $d^m \times 1$ and given in terms of the derivative vector $\mathbf{D}_{\mathbf{x}} = \left(\frac{\partial}{\partial x_1}, \frac{\partial}{\partial x_2}, \cdots, \frac{\partial}{\partial x_d}\right)'$  as
\begin{equation}
\mathbf{c}(m,d) = \left( \mathbf{D}_{\mathbf{x}}^{\otimes m} f(\mathbf{x})\right)\vline_{\mathbf{x}=\mathbf{0}}
\end{equation}
The Taylor series expansion near $\bs{\lambda} = \mathbf{0}$, called the Maclaurian series, of some required functions  based on the Equation  \eqref{ndTaylor} are derived in appendix \ref{Tayexpand}.  
\section{Moments, cumulants and characteristic function of a random vector}
\label{chfn}
Let $\mathbf{x} = \left(X_1, X_2,..., X_d\right)'$ be a d-dimensional random vector and $f(\mathbf{x})$ be its joint \textit{pdf} differentiable in each variable. 

\hspace{0.2 in}		The Characteristic function ($\mathcal{F}$) of $\mathbf{x}$ is defined as the expected value of $\mbox{e}^{i\mathbf{x}'\bs{\lambda}}$, where $\bs{\lambda} = (\lambda_1, \lambda_2, \ldots, \lambda_d)'$, $\bs{\lambda} \in \mathbb{R}^d$. 
Also, both the characteristic function and the \textit{pdf} are the Fourier Transform ($ \mathsf{F}$) of each other,  in the sense they are dual.
\begin{align}
\mathcal{F}_{\mathbf{x}}(\bs{\lambda})  = E\left\{\mbox{e}^{ i (\mathbf{x}'\bs{\lambda}) } \right\}  = \mathsf{F}(f(\mathbf{x})) =  \int_{-\infty}^{\infty}\cdots \int_{-\infty}^{\infty} f(\mathbf{x})\mbox{e}^{i(\mathbf{x}'\bs{\lambda})}d\mathbf{x}
\end{align}
Expanding $\mbox{e}^{i\mathbf{x}'\bs{\lambda}}$ using its Maclaurian series in  Equation \eqref{Tayeax} in appendix \ref{Tayexpand}, 
we get:
\begin{align}
\mathcal{F}_{\mathbf{x}}(\bs{\lambda}) = \sum_{k=0}^{\infty}  \mathbf{m}(k,d)'\frac{(i\bs{\lambda})^{\otimes k}}{k!}
\end{align} 
where, $\mathbf{m}(k,d)$ is the $k^{th}$ order moment vector of dimension $d^k \times 1$ and given by 
\begin{align}
\mathbf{m}(k,d) &= \int_{\mathbb{R}^{d}}\mathbf{x}^{\otimes k}f(\mathbf{x})d\mathbf{x} \nonumber \\ 
\label{deltafx}
\mbox{Also, } f(\mathbf{x}) = \mathsf{F}^{-1}(\mathcal{F}(\lambda)) &= \mathsf{F}^{-1}\left(\sum_{k=0}^{\infty}  \mathbf{m}(k,d)'\frac{(i\bs{\lambda})^{\otimes k}}{k!} \right) \\
 &= \sum_{k=0}^{\infty}  \frac{\mathbf{m}(k,d)'}{k!} \left( \frac{1}{(2\pi)^d} \int_{\mathbb{R}^{d}}(i\bs{\lambda)}^{\otimes k} \mbox{e}^{-i\mathbf{x}'\bs{\lambda}}d\bs{\lambda} \right) \nonumber \\
\label{fxdelta} 
 & = \sum_{k=0}^{\infty}  (-1)^k \frac{\mathbf{m}(k,d)'}{k!} \mathbf{D}^{(k)}\delta (\mathbf{x})  \\
 & \quad ( \because \mbox{Proof in Appendix \ref{prfderdel}}) \nonumber  
\end{align}
The Moment Generating Function (MGF) of $f(\mathbf{X})$ is given as
\begin{align}
\mathbf{M}(\bs{\lambda}) & = E\left\{\mbox{e}^{ \mathbf{x}'\bs{\lambda} } \right\} =  \int_{\mathbb{R}^{d}}f(\mathbf{X})\mbox{e}^{\mathbf{x}'\bs{\lambda}}d\mathbf{X}  \\
\label{ndMGF} 
 & = \sum_{k=0}^{\infty} \mathbf{m}(k,d)'\frac{\bs{\lambda}^{\otimes k}}{k!} \quad \mbox{   ( $\because$ Expanding $\mbox{e}^{\mathbf{x}'\bs{\lambda}}$)}
\end{align}
Assuming $\mathbf{M}(\bs{\lambda})$ and $\mathcal{F}(\bs{\lambda})$ are expanded using Taylor series,
\begin{equation}
\label{momkd}
\mathbf{m}(k,d) = \mathbf{D}_{\bs{\lambda}}^{\otimes k} \mathbf{M}(\bs{\lambda})\vline_{\bs{\lambda}=\mathbf{0}} = (-i)^k \mathbf{D}_{\bs{\lambda}}^{\otimes k} \mathcal{F}_{\mathbf{x}}(\bs{\lambda})\vline_{\bs{\lambda}=\mathbf{0}} 
\end{equation}
The Cumulant Generating Function (CGF) of $f(\mathbf{X})$ is given by,
\begin{equation}
\label{ndCGF}
\mathbf{C}(\bs{\lambda}) = \ln \mathbf{M}(\bs{\lambda}) = \sum_{k=1}^{\infty}\mathbf{c}(k,d)'\frac{\bs{\lambda}^{\otimes k}}{k!}
\end{equation}
where, $\mathbf{c}(k,d)$ is the $k^{th}$ order cumulant vector of dimension $d^k \times 1$. \\
The Cumulant Generating Function (CGF) of $f(\mathbf{X})$ can also be defined using the characteristic function, as under:
\begin{equation}
\label{ndCGF2}
\mathcal{C}(\bs{\lambda}) = \ln \mathcal{F}(\bs{\lambda}) = \sum_{k=1}^{\infty}\mathbf{c}(k,d)'\frac{(i\bs{\lambda})^{\otimes k}}{k!}
\end{equation}
Assuming $\mathbf{C}(\bs{\lambda})$ and $\mathcal{C}(\bs{\lambda})$ have been expanded using Taylor series, 
\begin{equation}
\label{cumkd}
\mathbf{c}(k,d) = \mathbf{D}_{\bs{\lambda}}^{\otimes k} \mathbf{C}(\bs{\lambda})\vline_{\bs{\lambda}=\mathbf{0}} = (-i)^k \mathbf{D}_{\bs{\lambda}}^{\otimes k} \mathcal{C}_{\mathbf{x}}(\bs{\lambda})\vline_{\bs{\lambda}=\mathbf{0}} 
\end{equation}
Significantly, this section has derived the moments and the cumulants of a random vector in vector notations. 
\subsection{Relation between the cumulant vectors and the moment vectors}
\label{RelCumMom}
The relation between the moments and the cumulants is given by combining Equation \eqref{ndMGF} and Equation \eqref{ndCGF} as below:
\begin{equation}
\label{momcumeq}
\mathbf{M}(\bs{\lambda}) = \sum_{k=0}^{\infty} \mathbf{m}(k,d)'\frac{\bs{\lambda}^{\otimes k}}{k!} = \exp \left( \sum_{k=1}^{\infty}\mathbf{c}(k,d)'\frac{\bs{\lambda}^{\otimes k}}{k!}\right)
\end{equation}
For $k=1$, using Equation \eqref{momkd}, we get: 
\begin{align*}
\mathbf{m}(1,d) = \mathbf{D}_{\bs{\lambda}}^{\otimes 1}  \mathbf{M}(\bs{\lambda})\vline_{\bs{\lambda}=\mathbf{0}} 	
\end{align*}
Applying K-derivative $(\mathbf{D}_{\bs{\lambda}}^{\otimes })$ to Equation \eqref{momcumeq}, 
\begin{align*}
\sum_{p=1}^{\infty} \mathbf{m}(p,d)'\frac{ \left( \bs{\lambda}^{\otimes (p-1)} \otimes \mathbf{I}_d \right)  }{(p-1)!}\vline_{\bs{\lambda} = \mathbf{0}} &= \sum_{p=1}^{\infty}\mathbf{c}(p,d)'\frac{ \left( \bs{\lambda}^{\otimes (p-1)} \otimes \mathbf{I}_d \right) }{(p-1)!} \nonumber \\
& \quad \otimes \exp \left( \sum_{q=1}^{\infty}\mathbf{c}(q,d)'\frac{\bs{\lambda}^{\otimes q}}{q!}\right)\vline_{\bs{\lambda}= \mathbf{0}} \\
\Rightarrow  \mathbf{m}(1,d) &= \mathbf{c}(1,d)
\end{align*}
Similarly, based on Equation \eqref{momkd}, for any k taking $k^{th}$ order K-derivative of above Equation  \eqref{momcumeq} on both sides relates $\mathbf{m}(k,d)$ with $\mathbf{c}(k,d)$.   For example, the cases for $k=2$ and  $k=3$ are shown in Appendix \ref{Apmomcum}. Overall, the vector moments in terms of the vector cumulants can be  summarized as under:
\begin{equation}
\label{eqndmom2cum}
\begin{aligned}
\mathbf{m}(1,d) &= \mathbf{c}(1,d)  \\
\mathbf{m}(2,d) &= \mathbf{c}(2,d) + \mathbf{c}(1,d)^{\otimes 2}  \\ 
\mathbf{m}(3,d) &= \mathbf{c}(3,d) + 3\mathbf{c}(2,d)\otimes \mathbf{c}(1,d) + \mathbf{c}(1,d)^{\otimes 3}  \\
\mathbf{m}(4,d) &=  \mathbf{c}(4,d) + 4\mathbf{c}(3,d)\otimes \mathbf{c}(1,d) + 3\mathbf{c}(2,d)^{\otimes 2} + 6 \mathbf{c}(2,d) \\
& \quad \otimes \mathbf{c}(1,d)^{\otimes 2}+ \mathbf{c}(1,d)^{\otimes 4}	 \\
\mathbf{m}(5,d) &= \mathbf{c}(5,d) + 5\mathbf{c}(4,d)\otimes \mathbf{c}(1,d)+ 10\mathbf{c}(3,d)\otimes \mathbf{c}(2,d) 
 \\
& \quad + 10\mathbf{c}(3,d) \otimes \mathbf{c}(1,d)^{\otimes 2} + 15\mathbf{c}(2,d)^{\otimes 2}\otimes \mathbf{c}(1,d)  \\
& + 10\mathbf{c}(2,d)\otimes \mathbf{c}(1,d)^{\otimes 3} + \mathbf{c}(1,d)^{\otimes 5}  \\
\mathbf{m}(6,d) &= \mathbf{c}(6,d) + 6\mathbf{c}(5,d)\otimes \mathbf{c}(1,d) + 15 \mathbf{c}(4,d) \otimes \mathbf{c}(2,d) \\
& \quad + 15 \mathbf{c}(4,d)  \otimes \mathbf{c}(1,d)^{\otimes 2}  + 10 \mathbf{c}(3,d)^{\otimes 2} + 60\mathbf{c}(3,d)\otimes \mathbf{c}(2,d) \\
& \quad \otimes \mathbf{c}(1,d) + 20\mathbf{c}(3,d)\otimes \mathbf{c}(1,d)^{\otimes 2} 
 + 15\mathbf{c}(2,d)^{\otimes 3}      \\
& + 45\mathbf{c}(2,d)^{\otimes 2}\otimes \mathbf{c}(1,d)^{\otimes 2} + 15\mathbf{c}(2,d)\otimes \mathbf{c}(1,d)^{\otimes 4} + \mathbf{c}(1,d)^{\otimes 6}   
\end{aligned}
\end{equation}
The above set of equations can be represented through more compact formulas as under. The Equation \eqref{ndmomcum} gives generalized $k^{th}$ order d-variate cumulant vector  in terms of the moment vectors and the 
Equation \eqref{ndcummom} gives the vice-a-versa. 
\begin{align}
\label{ndmomcum}
\mathbf{m}(k,d) =  \sum_{p=0}^{k-1} \binom{k-1}{p} \mathbf{K}^{-1}_{\mathfrak{p} h \leftrightarrow l}\mathbf{c}(k-p,d) \otimes \mathbf{m}(p,d) \\
\label{ndcummom}
\mathbf{c}(k,d) =  \mathbf{m}(k,d) - \sum_{p=1}^{k-1} \binom{k-1}{p} \mathbf{K}^{-1}_{\mathfrak{p} h \leftrightarrow l}\mathbf{c}(k-p,d) \otimes \mathbf{m}(p,d) 
\end{align} 
where, $\mathbf{K}^{-1}_{\mathfrak{p} h \leftrightarrow l}$ is a specific commutation matrix with corresponding dimensions that changes the place of the cumulants for Kronecker product such that the expression has decreasing order cumulants from left to the right, i.e. the higher order  cumulant vector on  the left and the lower order cumulant vector on the right. As the Kronecker products are non-commutative, without using the commutation matrices it would have been impossible to derive the compact formula. 
The derived  multivariate expressions reduce to the following expressions for dimension $d=1$ and are exactly same as those derived in   \citep{genGramC07}. 
\begin{equation}
\begin{aligned}
c(k,1) &=  m(k,1) - \sum_{p=1}^{k-1} \binom{k-1}{p} c(k-p,1)  m(p,1) \\ 
\mbox{or more simply, } c_k & =  m_k - \sum_{p=1}^{k-1} \binom{k-1}{p} c_{k-p} m_p 
\end{aligned}
\end{equation}
Thus, the  derived multivariate expressions in Equation \eqref{eqndmom2cum} are elementary vector extensions to those for univariate.  
\section{Multivariate \textit{pdf} representation in terms of the cumulants}
\label{pdfcum}
From Equation \eqref{deltafx} and Equation \eqref{ndCGF2}, the multivariate \textit{pdf} $f(\mathbf{x})$ can be written as:
\begin{align}
\label{ndpdf}
f(\mathbf{x}) = \mathsf{F}^{-1}(\mbox{e}^{\mathcal{C}(\bs{\lambda})}) = \left(\frac{1}{2\pi}\right)^d\int_{\mathbb{R}^{d}}\exp \left( \sum_{k=1}^{\infty}\mathbf{c}(k,d)'\frac{(i\bs{\lambda})^{\otimes k}}{k!} \right) \exp{(-i \mathbf{x}'\bs{\lambda})}d \bs{\lambda} 
\end{align}
As \textit{pdf} is a real function and $\mbox{Re}(e^{A+iB}e^{-iC}) = e^A\cos(B-C)= e^A\cos(C-B)$, the Equation \eqref{ndpdf} can be re-written as:
\begin{align} 
\label{ndpdfreal}
f(\mathbf{x}) &= \frac{1}{(2\pi)^d}\int_{\mathbb{R}^{d}} \exp{\left( \sum_{k=1}^{\infty}\frac{\mathbf{c}(2k,d)'}{2k!} {(i\bs{\lambda})^{\otimes 2k}} \right) } \nonumber \\
& \quad \cos\left( \mathbf{x}'\bs{\lambda} + \sum_{k=1}^{\infty}\frac{\mathbf{c}(2k-1,d)'}{(2k-1)!} (i)^{2k} (\bs{\lambda})^{\otimes (2k-1)} \right) d \bs{\lambda} 
\end{align}
The integrand in this equation is an even function. So, 
\begin{align}
\label{ndpdfeven}
f(\mathbf{x}) &= \frac{1}{(\pi)^d}\int_{(\mathbb{R}^{+})^d} \exp \left( \sum_{k=1}^{\infty}\frac{\mathbf{c}(2k,d)'}{2k!} {(i\bs{\lambda})^{\otimes 2k}} \right)  \nonumber \\
& \quad \cos\left( \mathbf{x}'\bs{\lambda} + \sum_{k=1}^{\infty}\frac{\mathbf{c}(2k-1,d)'}{(2k-1)!} (i)^{2k}(\bs{\lambda})^{\otimes (2k-1)}\right) d \bs{\lambda} 
\end{align} 
where, $\mathbb{R}^{+} =  \{ x \in \mathbb{R} : x \geq 0 \}$. 	
The Equation \eqref{ndpdfreal} and the Equation \eqref{ndpdfeven} give a multivariate \textit{pdf} in terms of the cumulants. As they are derived using Taylor series expansion, the infinite order differentiability is an implicit assumption. 

\hspace{0.2 in}	The equations can be verified using known \textit{pdf} examples with finite number of moments and cumulants.
 Let say, the 
 impulse delta density function has only the first order cumulant being non-zero and all other higher order cumulants are zero. Using this knowledge in Equation \eqref{ndpdfreal}, 
\begin{align}
f(\mathbf{x}) &= \frac{1}{(2\pi)^d}\int_{\mathbb{R}^{d}}  \cos\left( \left( \mathbf{x} - \mathbf{c}(1,d)\right)'\bs{\lambda} \right) d \bs{\lambda}  \\
&= \delta (\mathbf{x} - \mathbf{c}(1,d)) \quad \mbox{(shifted impulse delta function)} \nonumber 
\end{align}
Let's take another example, the Gaussian density function has first two order cumulants nonzero and all other order  cumulants are zero. Using this knowledge in Equation \eqref{ndpdfreal}, 
\begin{equation}
\label{Gpdf}
\begin{aligned}
f(\mathbf{x}) &= \frac{1}{(2\pi)^d}\int_{\mathbb{R}^{d}} \exp \left( - \frac{\mathbf{c}(2,d)'}{2}\bs{\lambda}^{\otimes 2} \right)  \cos\left( \left( \mathbf{x} - \mathbf{c}(1,d)\right)'\bs{\lambda} \right) d \bs{\lambda}  \\
&= G(\mathbf{x}) \quad  \mbox{ ( See Appendix \ref{Gaussint} for proof)}
\end{aligned}
\end{equation}
\section{The multivariate Hermite polynomials in integral form}
\label{DerGauss}
An interesting application of the integral form of multivariate \textit{pdf} representation is achieved in this section. 
The multivariate Gaussian expressed as in Equation \eqref{Gpdf} is used to derive it's differentials and Hermite  polynomials in a simple way. Taking $k^{th}$ K-derivative of $G(\mathbf{x})$, 
\begin{align}
\label{derkGx}
G^{(k)}(\mathbf{x}) 
& := \mathbf{D}_{\mathbf{x}}^{\otimes k}G(\mathbf{x}) 
= \frac{1}{(2\pi)^d}\int_{\mathbb{R}^{d}} \bs{\lambda}^{\otimes k} \exp \left( -\frac{\mathbf{c}(2,d)'}{2}\bs{\lambda}^{\otimes 2} \right)  \nonumber \\
& \quad  \cos\left( \left( \mathbf{x} - \mathbf{c}(1,d) \right)'\bs{\lambda} + \frac{k\pi}{2} \right) d \bs{\lambda} 
\end{align}
The multivariate Hermite polynomials defined by \citet{dnHermite96} are defined as under: 
\begin{align}
\label{hermitend}
\mathbf{H}_k(\mathbf{x}; \mathbf{0},\mathbf{C}_{\mathbf{x}}) &= [G(\mathbf{x};\mathbf{0},\mathbf{C}_{\mathbf{x}})]^{-1} (-1)^{k} \left( \mathbf{C}_{\mathbf{x}} \mathbf{D}_x \right)^{\otimes k} G(\mathbf{x};\mathbf{0},\mathbf{C}_{\mathbf{x}}) \\
\mbox{where, } G(\mathbf{x};\mathbf{0},\mathbf{C}_{\mathbf{x}}) &= \left|\mathbf{C}_{\mathbf{x}} \right|^{-1/2}(2\pi)^{-d/2}\exp\left(- \frac{1}{2}\mathbf{x}' \mathbf{C}_{\mathbf{x}}^{-1} \mathbf{x} \right) \nonumber \\  & = \left|\mathbf{C}_{\mathbf{x}} \right|^{-1/2}(2\pi)^{-d/2}\exp\left( - \frac{1}{2} \left( Vec \mbox{ } \mathbf{C}_{\mathbf{x}}^{-1}\right)' \mathbf{x}^{\otimes 2} \right)
\end{align}
This is equivalent to the 1-dimensional definition of Hermite polynomials\footnote{This is the 'probabilists' Hermite polynomials and not the 'physicists' Hermite polynomials used by \citet{genGramC07}.}  by Rodrigues's formula in Equation \eqref{hermite1d}, except the introduction of matrix $\mathbf{C}_{\mathbf{x}}$.
\begin{align}
\label{hermite1d}
H_k(x) = [G(x)]^{-1}(-1)^k \frac{d^k}{dx^k}G(x) \mbox{ where, }  G(x) = \frac{1}{2\pi}\mbox{e}^{-\frac{1}{2}x^2} 
\end{align}
Using Equation \eqref{derkGx}, the Equation \eqref{hermitend} for multivariate Hermite polynomials is rewritten as: 
\begin{align}
\label{hermitenew}
\mathbf{H}_k & (\mathbf{x}; \mathbf{0},\mathbf{C}_{\mathbf{x}} )
 = (2\pi)^{d/2} \left|\mathbf{C}_{\mathbf{x}}\right|^{1/2} (-1)^k \left(\mathbf{C}_{\mathbf{x}} \right)^{\otimes k} \exp \left( \frac{1}{2} \left(Vec \mbox{ } \mathbf{C}_{\mathbf{x}}^{-1}\right)' \mathbf{x}^{\otimes 2} \right) \frac{1}{(2\pi)^d} \nonumber \\
&\quad  \int_{\mathbb{R}^{d}} \bs{\lambda}^{\otimes k} \exp \left( -\frac{\mathbf{c}(2,d)'}{2}\bs{\lambda}^{\otimes 2} \right)  \cos\left( (\mathbf{x} - \mathbf{c}(1,d))'\bs{\lambda} + \frac{k\pi}{2} \right) d \bs{\lambda} 
\end{align}
Taking  $\mathbf{C}_{\mathbf{x}} = \mathbf{I}_d$, where $\mathbf{I}_d$ is $d \times d$  identity matrix; using the property $(Vec \mbox{ } \mathbf{I}_d )'\mathbf{x}^{\otimes 2} = \mathbf{x}'\mathbf{x}$ and using change of variable as $\bs{\lambda}/\sqrt{2} =  \mathbf{u}$ -   the integral form of multivariate Hermite polynomials is obtained as under:
\begin{align}
\label{hermitend01}
\mathbf{H}_k(\mathbf{x}; \mathbf{0},\mathbf{I}_d) & = (2)^{ \frac{k+d+1}{2} } (\pi)^{-d/2}
\exp \left( \frac{1}{2}\mathbf{x}' \mathbf{x}\right) \nonumber\\
& \quad \int_{\mathbb{R}^{d}} \mathbf{u}^{\otimes k} 
\exp \left( -\mathbf{u}'\mathbf{u} \right)  \cos\left( \sqrt{2}\mathbf{x}'\mathbf{u} - \frac{k\pi}{2} \right) d \mathbf{u} 
\end{align}
    
\hspace{0.2 in}	 The result in Equation \eqref{derkGx} can also be obtained  using Equation \eqref{Gpdf} and applying the derivative property of Fourier transform ($\mathsf{F}$). The $k^{th}$ derivative of $G(x)$ is given by, 
\begin{align}
\label{eqFgk}
G^{(k)}(\mathbf{x}) &= \mathsf{F}^{-1}\left( \mathsf{F}(G^{(k)}(\mathbf{x}))\right) \\
&= \frac{1}{(2\pi)^d} \int_{\mathbb{R}^{d}} \left( i\bs{\lambda} \right)^{\otimes k}\mathsf{F}(G(\mathbf{x})) \exp^{-i\mathbf{x}'\bs{\lambda}} d \bs{\lambda} \\
&= \frac{1}{(2\pi)^d}\int_{\mathbb{R}^{d}} \bs{\lambda}^{\otimes k} \exp \left( -\frac{\mathbf{c}(2,d)'}{2}\bs{\lambda}^{\otimes 2} \right)  \cos\left( \left( \mathbf{x} - \mathbf{c}(1,d) \right)'\bs{\lambda} + \frac{k\pi}{2} \right) d \bs{\lambda} 
\end{align}
\section[Multivariate Gram-Charlier A series ]{Multivariate Gram-Charlier A series} 
\label{ndGCAseries}
Till now, the article has derived - an unknown \textit{pdf} expressed in terms of its cumulants  in Equation \eqref{ndpdfreal};  the Gaussian density function expressed in terms of its cumulants in Equation \eqref{Gpdf} and the Hermite polynomials  in Equation \eqref{hermitend01}. Based on them, the multivariate Gram Charlier A series that expresses an unknown \textit{pdf} using Gaussian density as a reference can be obtained. The expansion assumes  first and second order cumulants being same for both the unknown \textit{pdf} and the reference \textit{pdf}. Using the expansion $\exp(A+B)\cos(C+D) = \exp(A)\exp(B)(\cos C \cos D -\sin C \sin D)$, the Equation \eqref{ndpdfreal} can be re-written as:
\begin{align}
\label{ndpdfGx}
f(\mathbf{x}) &= \frac{1}{(2\pi)^d}\int_{\mathbb{R}^{d}} 
\exp \left(- \frac{\mathbf{c}(2,d)'}{2}\bs{\lambda}^{\otimes 2} \right) 
\exp{ \left( \sum_{k=2}^{\infty}\frac{\mathbf{c}(2k,d)'}{2k!} {(i\bs{\lambda})^{\otimes 2k}} \right) } \nonumber \\ 
& \quad \left\{  
\cos( ( \mathbf{x} - \mathbf{c}(1,d) )'\bs{\lambda} ) \cos \left( \sum_{k=2}^{\infty} \frac{ \mathbf{c}(2k-1,d)' }{(2k-1)!} (i)^{2k} (\bs{\lambda})^{\otimes (2k-1)} \right) \right. \nonumber \\
& \quad \left. -  \sin( (\mathbf{x} - \mathbf{c}(1,d))'\bs{\lambda}) \sin \left( \sum_{k=2}^{\infty}\frac{\mathbf{c}(2k-1,d)'}{(2k-1)!} (i)^{2k} (\bs{\lambda})^{\otimes (2k-1)} \right) \right\} d \bs{\lambda} 
\end{align}   
Using the expansions in Appendix \ref{Tayexpand}, parts of the Equation \eqref{ndpdfGx} can be approximated upto maximum $6^{th}$-order statistics 
as under: 
\begin{align}
\label{excoscum}
 \exp & \left( \sum_{k=2}^{\infty}\frac{\mathbf{c}(2k,d)'}{2k!} (i\bs{\lambda})^{\otimes 2k} \right) \cos \left( \sum_{k=2}^{\infty} \frac{ \mathbf{c}(2k-1,d)' }{(2k-1)!} (i)^{2k} (\bs{\lambda})^{\otimes (2k-1)} \right) \nonumber \\
& \quad =  1 + \left( \frac{\mathbf{c}(4,d)'\bs{\lambda}^{\otimes 4}}{4!} -  \frac{\mathbf{c}(6,d)'\bs{\lambda}^{\otimes 6}}{6!}\right) - \frac{1}{2}\left( \frac{\mathbf{c}(3,d)'\bs{\lambda}^{\otimes 3}}{3!} \right)^{\otimes 2} + \ldots \\
\label{exsincum}
 \exp & \left( \sum_{k=2}^{\infty} \frac{\mathbf{c}(2k,d)'}{2k!} (i\bs{\lambda})^{\otimes 2k} \right) \sin \left( \sum_{k=2}^{\infty}\frac{\mathbf{c}(2k-1,d)'}{(2k-1)!} (-i)^{2k} (\bs{\lambda})^{\otimes (2k-1)} \right) \nonumber \\
&\quad = \frac{\mathbf{c}(3,d)'\bs{\lambda}^{\otimes 3}}{3!} -  \frac{\mathbf{c}(5,d)'\bs{\lambda}^{\otimes 5}}{5!} + \ldots 
\end{align}
Using above Equation \eqref{excoscum} and Equation \eqref{exsincum}, the Equation \eqref{ndpdfGx} can be re-written as:
\begin{align}
f(\mathbf{x}) &= \frac{1}{(2\pi)^d} \int_{\mathbb{R}^{d}} \left\{
\exp \left( - \frac{\mathbf{c}(2,d)'}{2}\bs{\lambda}^{\otimes 2} \right) \cos( ( \mathbf{x} - \mathbf{c}(1,d) )'\bs{\lambda} ) \right. \nonumber \\
& \quad \left. \left(   1 + \frac{\mathbf{c}(4,d)'\bs{\lambda}^{\otimes 4}}{4!} - \frac{1}{6!} \left(  \mathbf{c}(6,d) - 10 \mathbf{c}(3,d)^{\otimes 2} \right)'{\bs{\lambda}^{\otimes 6} } + \ldots \right) \right\} d \bs{\lambda} \nonumber \\
& \quad -  \frac{1}{(2\pi)^d}\int_{\mathbb{R}^{d}} \left\{
\exp \left(- \frac{\mathbf{c}(2,d)'}{2}\bs{\lambda}^{\otimes 2} \right) \cos \left( \left(\mathbf{x} - \mathbf{c}(1,d)\right)'\bs{\lambda} - \frac{\pi}{2}\right) \right. \nonumber \\
& \quad \left. \left(  \frac{\mathbf{c}(3,d)'\bs{\lambda}^{\otimes 3}}{3!} -  \frac{\mathbf{c}(5,d)'\bs{\lambda}^{\otimes 5}}{5!} + \ldots 
 \right) \right\} d \bs{\lambda} 
\end{align}
Using the Equation \eqref{derkGx} for derivatives of Gaussian defined, the above Equation can be simplified as:
\begin{align} 
\label{GCA}
f(\mathbf{x}) & = G(\mathbf{x}) -   \frac{\mathbf{c}(3,d)'}{3!} G^{(3)}(\mathbf{x}) +  \frac{\mathbf{c}(4,d)'}{4!}G^{(4)}(\mathbf{x}) - \frac{\mathbf{c}(5,d)'}{5!} G^{(5)}(\mathbf{x}) \nonumber \\
& \quad + \frac{\mathbf{c}(6,d)' + 10 \mathbf{c}(3,d)^{\otimes 2'}}{6!} G^{(6)}(\mathbf{x}) + \ldots
\end{align}
The Equation \eqref{GCA} is the Gram-Charlier A series expressed directly in terms of the cumulants and the derivatives  of the Gaussian \textit{pdf}. Usually, the GCA is represented in terms of the Hermite polynomials.  So, the GCA expansion (Equation \eqref{GCA}) in terms of the  Hermite polynomials; either using definition in Equation \eqref{hermitend} or using $\mathbf{H}_k( \mathbf{x}; \mathbf{0}, \mathbf{C}_{\mathbf{x}}^{-1})$ derived in Equation \eqref{hermitenew}; can be re-written as: 
\begin{align} 
\label{GCAhermitecx}
f(\mathbf{x}) & = G(\mathbf{x}) \left[ 1 + \frac{\mathbf{c}(3,d)'}{3!} \left( \mathbf{C}_{\mathbf{x}}^{-1} \right)^{\otimes 3} \mathbf{H}_3(\mathbf{x}; \mathbf{0},\mathbf{C}_{\mathbf{x}} ) +  \frac{\mathbf{c}(4,d)'}{4!}\left( \mathbf{C}_{\mathbf{x}}^{-1} \right)^{\otimes 4}\mathbf{H}_4(\mathbf{x}; \mathbf{0},\mathbf{C}_{\mathbf{x}} )   \right. \nonumber \\
& \quad \left. + \frac{\mathbf{c}(5,d)'}{5!} \left( \mathbf{C}_{\mathbf{x}}^{-1} \right)^{\otimes 5}\mathbf{H}_5(\mathbf{x}; \mathbf{0},\mathbf{C}_{\mathbf{x}} ) 
  + \frac{\mathbf{c}(6,d)' + 10 \mathbf{c}(3,d)^{\otimes 2'}}{6!} \left( \mathbf{C}_{\mathbf{x}}^{-1} \right)^{\otimes 6} \right. \nonumber \\
& \quad \left. \mathbf{H}_6(\mathbf{x}; \mathbf{0},\mathbf{C}_{\mathbf{x}} ) + \ldots  \right] \\
\label{GCAhermite}
f(\mathbf{x}) & = G(\mathbf{x}) \left[ 1 + \frac{\mathbf{c}(3,d)'}{3!} \mathbf{H}_3(\mathbf{x}; \mathbf{0},\mathbf{I}_d ) +  \frac{\mathbf{c}(4,d)'}{4!}\mathbf{H}_4(\mathbf{x}; \mathbf{0},\mathbf{I}_d ) + \frac{\mathbf{c}(5,d)'}{5!} \mathbf{H}_5(\mathbf{x}; \mathbf{0},\mathbf{I}_d ) \right. \nonumber \\
& \quad \left.  + \frac{\mathbf{c}(6,d)' + 10 \mathbf{c}(3,d)^{\otimes 2'}}{6!} \mathbf{H}_6(\mathbf{x}; \mathbf{0},\mathbf{I}_d ) + \ldots  \right] 
\end{align}
Finally, the GCA series, in vector notations, can be expressed either using $k^{th}$ order derivative of Gaussian $(G_k(\mathbf{x}))$ or using $k^{th}$ order  vector Hermite polynomials $(\mathbf{H}_k(\mathbf{x}))$ as under:  
\begin{align}
\label{ndGCAshort}
f(\mathbf{x}) &= \sum_{k=0}^{\infty} (-1)^k\frac{\mathbf{c}(k,d)'}{k!}G^{(k)}(\mathbf{x}) \\
\label{ndGCAhershort}
 & = \sum_{k=0}^{\infty} \frac{\mathbf{c}(k,d)'}{k!}\mathbf{H}_k(\mathbf{x}; \mathbf{0},\mathbf{I}_d)
\end{align}
\section{Multivariate Generalized Gram-Charlier series}
\label{ndGGC}
To derive the generalized Gram-Charlier series, an unknown \textit{pdf} $f(\mathbf{x})$ need be represented in terms of any known  reference \textit{pdf} $\psi(\mathbf{x})$, where both the \textit{pdf}s are represented in terms of their cumulants.  
Let the $k^{th}$ order cumulant vector of the reference \textit{pdf} $\psi(\mathbf{x})$ be $\mathbf{c}_r(k,d)$. Then, the $k^{th}$ order cumulant difference vector $\bs{\delta}(k,d)$ is:  $\bs{\delta}(k,d) = \mathbf{c}(k,d) - \mathbf{c}_r(k,d), \forall k$. Using $\bs{\delta}(k,d)$, the Equation \eqref{ndpdfreal} can be re-written 
as under:
\begin{align}
f(\mathbf{x}) &= \frac{1}{(2\pi)^d}\int_{\mathbb{R}^{d}} \exp{\left( \sum_{k=1}^{\infty}\frac{\mathbf{c}_r(2k,d)'}{2k!} {(i\bs{\lambda})^{\otimes 2k}} \right) }\exp{\left( \sum_{k=1}^{\infty}\frac{\bs{\delta}(2k,d)'}{2k!} {(i\bs{\lambda})^{\otimes 2k}} \right) } \nonumber \\
 & \quad \cos\left( \mathbf{x}'\bs{\lambda} + \left( \sum_{k=1}^{\infty}\frac{\mathbf{c}_r(2k-1,d)}{(2k-1)!}  + \sum_{k=1}^{\infty}\frac{\bs{\delta}(2k-1,d)}{(2k-1)!} \right)' (i)^{2k} (\bs{\lambda})^{\otimes (2k-1)} \right) d \bs{\lambda} \\
\label{ndpdfcrdelta} 
& =\frac{1}{(2\pi)^d}\int_{\mathbb{R}^{d}}  \exp{\left( \sum_{k=1}^{\infty}\frac{\mathbf{c}_r(2k,d)'}{2k!} {(i\bs{\lambda})^{\otimes 2k}} \right) } \exp{\left( \sum_{k=1}^{\infty}\frac{\bs{\delta}(2k,d)'}{2k!} {(i\bs{\lambda})^{\otimes 2k}} \right) }  \nonumber \\ & \quad 
\left\{ \cos\left( \mathbf{x}'\bs{\lambda} + 
 \sum_{k=1}^{\infty}\frac{\mathbf{c}_r(2k-1,d)'}{(2k-1)!} (i)^{2k} (\bs{\lambda})^{\otimes (2k-1)} \right)
\cos\left( \sum_{k=1}^{\infty}\frac{\bs{\delta}(2k-1,d)'}{(2k-1)!} \right. \right. \nonumber \\
& \left. \left. (i)^{2k} (\bs{\lambda})^{\otimes (2k-1)} \right) 
- \sin\left( \mathbf{x}'\bs{\lambda} + 
 \sum_{k=1}^{\infty}\frac{\mathbf{c}_r(2k-1,d)'}{(2k-1)!} (i)^{2k} (\bs{\lambda})^{\otimes (2k-1)} \right)
\right. \nonumber \\  
& \quad \left. \sin\left( \sum_{k=1}^{\infty}\frac{\bs{\delta}(2k-1,d)'}{(2k-1)!} 
(i)^{2k} (\bs{\lambda})^{\otimes (2k-1)} \right) \right\}d\bs{\lambda}
\end{align}
Using the expansions in Appendix \ref{Tayexpand}, parts of the Equation \eqref{ndpdfcrdelta} can be approximated upto maximum $6^{th}$-order statistics 
as under: 
\begin{align}
\label{expcos}
\exp & {\left( \sum_{k=1}^{\infty}\frac{\bs{\delta}(2k,d)'}{2k!} {(i\bs{\lambda})^{\otimes 2k}} \right) }\cos\left( \sum_{k=1}^{\infty}\frac{\bs{\delta}(2k-1,d)'}{(2k-1)!} (i)^{2k} (\bs{\lambda})^{\otimes (2k-1)} \right)  \nonumber \\
 &=   1 + \left( \sum_{k=1}^{\infty}\frac{\bs{\delta}(2k,d)'}{2k!} {(i\bs{\lambda})^{\otimes 2k}} \right)+ \frac{1}{2}\left( \sum_{k=1}^{\infty}\frac{\bs{\delta}(2k,d)'}{2k!} {(i\bs{\lambda})^{\otimes 2k}} \right)^{\otimes 2} \nonumber \\ 
 & \quad  - \frac{1}{2}
\left( \sum_{k=1}^{\infty}\frac{\bs{\delta}(2k-1,d)'}{(2k-1)!} (i)^{2k} (\bs{\lambda})^{\otimes (2k-1)} \right)^{\otimes 2} - \ldots \nonumber \\
& =  1 + \left( - \frac{\bs{\delta}(2,d)'}{2!}\bs{\lambda}^{\otimes 2} + \frac{\bs{\delta}(4,d)'}{4!}\bs{\lambda}^{\otimes 4} -  \frac{\bs{\delta}(6,d)'}{6!}\bs{\lambda}^{\otimes 6} \right) + \frac{1}{2} \left( - \frac{\bs{\delta}(2,d)'}{2!}\bs{\lambda}^{\otimes 2}  \right)^{\otimes 2} 
\nonumber \\ 
&\quad  - \frac{1}{2}\left( - \bs{\delta}(1,d)'\bs{\lambda} 
+ \frac{\bs{\delta}(3,d)'}{3!}\bs{\lambda}^{\otimes 3}  \right)^{\otimes 2}  + \frac{1}{6}\left( - \frac{\bs{\delta}(2,d)'}{2!}\bs{\lambda}^{\otimes 2} \right)^{\otimes 3} 
- \frac{1}{2}\left( \left( - \frac{\bs{\delta}(2,d)'}{2!}\bs{\lambda}^{\otimes 2} \right. \right.
\nonumber \\ 
& \quad \left. \left. + \frac{\bs{\delta}(4,d)'}{4!}\bs{\lambda}^{\otimes 4} \right)\otimes  \left( \bs{\delta}(1,d)'\bs{\lambda}\right)^{\otimes 2} \right) 
- \frac{1}{4}\left( \left(- \frac{\bs{\delta}(2,d)'}{2!}\bs{\lambda}^{\otimes 2}\right)^{\otimes 2} \otimes  \left( \bs{\delta}(1,d)'\bs{\lambda}\right)^{\otimes 2}\right) \nonumber \\
& \quad +  \frac{1}{4!} \left( - \frac{\bs{\delta}(2,d)'}{2!}\bs{\lambda}^{\otimes 2}  \right)^{\otimes 4} 
- \frac{1}{4!}\left( - \bs{\delta}(1,d)'\bs{\lambda} \right)^{\otimes 4} -  \left( \frac{\bs{\delta}(2,d)'}{2!}\bs{\lambda}^{\otimes 2} \otimes \bs{\delta}(1,d)'\bs{\lambda}^{\otimes 4} \right) + 
\ldots \nonumber \\
 & =   1 - \frac{1}{2}\left( \bs{\delta}(1,d)^{\otimes 2} + \bs{\delta}(2,d)\right)'\bs{\lambda}^{\otimes 2} 
 + \frac{1}{4!} \left( \bs{\delta}(1,d)^{\otimes 4}+ 6 \bs{\delta}(2,d)\otimes \bs{\delta}(1,d)^{\otimes 2}  \right. \nonumber \\  & \quad  \left. + 3\bs{\delta}(2,d)^{\otimes 2}
+ 4\bs{\delta}(3,d)\otimes \bs{\delta}(1,d) + \bs{\delta}(4,d)\right)'\bs{\lambda}^{\otimes 4} 
- \frac{1}{6!}\left( \bs{\delta}(1,d)^{\otimes 6} + 15\bs{\delta}(2,d)^{\otimes 3} \right. \nonumber \\ & \quad \left. 
+ 10 \bs{\delta}(3,d)^{\otimes 2} + 
15 \bs{\delta}(4,d) \otimes \bs{\delta}(2,d) 
+ 15 \bs{\delta}(4,d) 
\otimes \bs{\delta}(1,d)^{\otimes 2}  
  + 20 \bs{\delta}(3,d)\otimes \bs{\delta}(1,d)^{\otimes 2} \right. \nonumber \\ &\quad \left.
+ 15\bs{\delta}(2,d)\otimes \bs{\delta}(1,d)^{\otimes 4} + 45\bs{\delta}(2,d)^{\otimes 2}\otimes \bs{\delta}(1,d)^{\otimes 2} 
+ 6 \bs{\delta}(5,d)\otimes \bs{\delta}(1,d) \right. \nonumber \\ &\quad \left.
+ 60\bs{\delta}(3,d)\otimes \bs{\delta}(2,d)\otimes \bs{\delta}(1,d) + \bs{\delta}(6,d)\right)'\bs{\lambda}^{\otimes 6} + \ldots 
\end{align}
Similarly, 
\begin{align}
\label{expsin}
 \exp & {\left( \sum_{k=1}^{\infty}\frac{\bs{\delta}(2k,d)'}{2k!} {(i\bs{\lambda})^{\otimes 2k}} \right) }\sin\left( \sum_{k=1}^{\infty}\frac{\bs{\delta}(2k-1,d)'}{(2k-1)!} (-i)^{2k} (\bs{\lambda})^{\otimes (2k-1)} \right)  \nonumber \\
 &=  \left( - \bs{\delta}(1,d)'\bs{\lambda} + \frac{ \bs{\delta}(3,d)'}{3!}\bs{\lambda}^{\otimes 3} - \frac{ \bs{\delta}(5,d)'}{5!}\bs{\lambda}^{\otimes 5} \right)
+ \left( \frac{ \bs{\delta}(2,d)'\bs{\lambda}^{\otimes 2}}{2!} \otimes \frac{ \bs{\delta}(1,d)'\bs{\lambda}}{1} \right. \nonumber \\
& \quad \left.  + \frac{ \bs{\delta}(3,d)'\bs{\lambda}^{\otimes 3}}{3!}\otimes \frac{ \bs{\delta}(2,d)'\bs{\lambda}^{\otimes 2}}{2!} + \frac{ \bs{\delta}(4,d)'\bs{\lambda}^{\otimes 4}}{4!}\otimes \frac{ \bs{\delta}(1,d)'\bs{\lambda}^{\otimes 1}}{1!} \right) \nonumber \\
& \quad - \frac{1}{2}\left( \frac{ \bs{\delta}(2,d)'\bs{\lambda}^{\otimes 2} }{2!} \right)^{\otimes 2}\otimes \frac{ \bs{\delta}(1,d)'\bs{\lambda}^{\otimes 1}}{1!} - \left( - \bs{\delta}(1,d)'\bs{\lambda}\right)^{\otimes 3} +  \ldots\nonumber \\
& =  \bs{\delta}(1,d)'\bs{\lambda} + \frac{1 }{3!}\left( \bs{\delta}(3,d) + 3\bs{\delta}(2,d)\otimes \bs{\delta}(1,d) + \bs{\delta}(1,d)^{\otimes 3} \right)' \bs{\lambda}^{\otimes 3} 
- \frac{ 1}{5!}\left( \bs{\delta}(5,d) \right. \nonumber \\ &\quad \left.
+ 5\bs{\delta}(4,d)\otimes \bs{\delta}(1,d)+ 15\bs{\delta}(2,d)^{\otimes 2}\otimes \bs{\delta}(1,d) + 10\bs{\delta}(3,d)\otimes \bs{\delta}(2,d) \right. \nonumber \\
& \quad \left. + 10 \bs{\delta}(2,d)\otimes \bs{\delta}(1,d)^{\otimes 3} + 10 \bs{\delta}(3,d)\otimes \bs{\delta}(1,d)^{\otimes 2} + \bs{\delta}(1,d)^{\otimes 5} \right)'\bs{\lambda}^{\otimes 5} \ldots 
\end{align}
Now, $\mathbf{D}_{\mathbf{x}}^{\otimes k} \psi(\mathbf{x})$ can be obtained by taking $k^{th}$-order K-derivative of Equation \eqref{ndpdfreal} as under: 
\begin{align} 
\label{derkref}
\psi^{(k)}(\mathbf{x}) & = \mathbf{D}_{\mathbf{x}}^{\otimes k} \psi(\mathbf{x}) = \frac{1}{(2\pi)^d}\int_{\mathbb{R}^{d}} (\bs{\lambda})^{\otimes k} \exp{\left( \sum_{k=1}^{\infty}\frac{\mathbf{c}_r(2k,d)'}{2k!} {(i\bs{\lambda})^{\otimes 2k}} \right) }  \nonumber \\
& \quad \quad \cos\left( \mathbf{x}'\bs{\lambda} + \sum_{k=1}^{\infty}\frac{\mathbf{c}_r(2k-1,d)'}{(2k-1)!} (i)^{2k} (\bs{\lambda})^{\otimes (2k-1)} + \frac{k\pi}{2} \right) d \bs{\lambda} 
\end{align}
The above Equation \eqref{derkref} with the previous results on expansions in Equation \eqref{expcos} and Equation \eqref{expsin} can be used to simplify the Equation \eqref{ndpdfcrdelta}. 
This derives the Generalized Gram-Charlier (GGC) series expressing an unknown \textit{pdf} $f(\mathbf{x})$ in terms of the cumulant difference vectors $(\bs{\delta}(k,d))$ and derivatives of a reference \textit{pdf} $\psi^{(k)}(\mathbf{x})$ as under: 
\begin{align}
\label{GGC}
f(\mathbf{x}) = \sum_{k=0}^{\infty} (-1)^k\frac{\bs{\alpha}(k,d)'}{k!}\psi^{(k)}(\mathbf{x})
\end{align}
where,
\begin{equation}
\label{delta2alpha}
\begin{aligned}
\bs{\alpha}(0,d) &= 1  \\
\bs{\alpha}(1,d) &= \bs{\delta}(1,d)  \\
\bs{\alpha}(2,d) &= \bs{\delta}(2,d) + \bs{\delta}(1,d)^{\otimes 2}  \\ 
\bs{\alpha}(3,d) &= \bs{\delta}(3,d) + 3\bs{\delta}(2,d)\otimes \bs{\delta}(1,d) + \bs{\delta}(1,d)^{\otimes 3}  \\
\bs{\alpha}(4,d) &=  \bs{\delta}(4,d) + 4\bs{\delta}(3,d)\otimes \bs{\delta}(1,d) + 3\bs{\delta}(2,d)^{\otimes 2} + 6 \bs{\delta}(2,d) \\
& \otimes \bs{\delta}(1,d)^{\otimes 2} + \bs{\delta}(1,d)^{\otimes 4}	 \\
\bs{\alpha}(5,d) &= \bs{\delta}(5,d) + 5\bs{\delta}(4,d)\otimes \bs{\delta}(1,d)+ 10\bs{\delta}(3,d)\otimes \bs{\delta}(2,d) 	\\
& \quad + 10\bs{\delta}(3,d) \otimes \bs{\delta}(1,d)^{\otimes 2}  
  + 15\bs{\delta}(2,d)^{\otimes 2}\otimes \bs{\delta}(1,d)    \\
  & + 10\bs{\delta}(2,d)\otimes \bs{\delta}(1,d)^{\otimes 3}+ \bs{\delta}(1,d)^{\otimes 5}  \\
\bs{\alpha}(6,d) &= \bs{\delta}(6,d) + 6\bs{\delta}(5,d)\otimes \bs{\delta}(1,d) + 15 \bs{\delta}(4,d) \otimes \bs{\delta}(2,d)  \\
& \quad + 15 \bs{\delta}(4,d) \otimes \bs{\delta}(1,d)^{\otimes 2}  + 10 \bs{\delta}(3,d)^{\otimes 2}  + 60\bs{\delta}(3,d)\otimes \bs{\delta}(2,d) \\
& \quad \otimes \bs{\delta}(1,d) + 20\bs{\delta}(3,d) \otimes \bs{\delta}(1,d)^{\otimes 2} + 15\bs{\delta}(2,d)^{\otimes 3}  + 45\bs{\delta}(2,d)^{\otimes 2} \\
& \quad \otimes \bs{\delta}(1,d)^{\otimes 2} + 15\bs{\delta}(2,d)\otimes \bs{\delta}(1,d)^{\otimes 4} + \bs{\delta}(1,d)^{\otimes 6}   
\end{aligned}
\end{equation}
The above set of equations \eqref{delta2alpha} has exact resemblance with that expressing moments in terms of the cumulants in Section \ref{RelCumMom}. 
This must happen, as Equation \eqref{GGC} for GGC expansion with $\delta(\mathbf{x})$ as a reference \textit{pdf}
is matching Equation \eqref{fxdelta}.  
This matching proves that $\bs{\alpha}(k,d)$  is related in same way to $\bs{\delta}(k,d)$, as $\mathbf{m}(k,d)$ to $\mathbf{c}(k,d)$. That is,:
\begin{align}
\label{alphadeltaeq}
\sum_{k=0}^{\infty} \bs{\alpha}(k,d)'\frac{\bs{\lambda}^{\otimes k}}{k!} = \exp \left( \sum_{k=1}^{\infty}\bs{\delta}(k,d)'\frac{\bs{\lambda}^{\otimes k}}{k!}\right)
\end{align}
Further, the $\bs{\alpha}(k,d)$ in Equation \eqref{GGC} recursively can be obtained in terms of the cumulant difference vector $\bs{\delta}(k,d)$ as  under:
\begin{align}
\label{ndalphadelta}
\bs{\alpha}(k,d) =  \sum_{p=0}^{k-1} \binom{k-1}{p} \mathbf{K}^{-1}_{\mathfrak{p} h \leftrightarrow l}\bs{\delta}(k-p,d) \otimes \bs{\alpha}(p,d) 
\end{align}
where, $\mathbf{K}^{-1}_{\mathfrak{p} h \leftrightarrow l}$ is a specific commutation matrix; as described previously; to    
change the order of the cumulants for Kronecker product such that the expression has decreasing order cumulants from left to the right. 

\hspace{0.2 in}		The verification of the derived GGC can be obtained by taking Gaussian density as a reference \textit{pdf}.  With Gaussian density as  a reference, $\bs{\delta}(1,d)= \mathbf{0}$ and $\bs{\delta}(2,d) = \mathbf{0}$. So, the coefficients $\bs{\alpha}(k,d)$ in Equation \eqref{GGC} can be derived as under:
\begin{equation}
\begin{aligned}
\bs{\alpha}(0,d) &= 1 \\
\bs{\alpha}(1,d) &= \mathbf{0} \\
\bs{\alpha}(2,d) &= \mathbf{0} \\ 
\bs{\alpha}(3,d) &= \bs{\delta}(3,d) = \mathbf{c}(3,d)  \\
\bs{\alpha}(4,d) &=  \bs{\delta}(4,d)  = \mathbf{c}(4,d) 	\\
\bs{\alpha}(5,d) &= \bs{\delta}(5,d) =  \mathbf{c}(5,d) \\
\bs{\alpha}(6,d) &= \bs{\delta}(6,d) + 10 \bs{\delta}(3,d)^{\otimes 2}   = \mathbf{c}(6,d) + 10 \mathbf{c}(3,d)^{\otimes 2}
\end{aligned}
\end{equation}
Thus, the GGC series is derived and verified using known examples. 
\section{Characteristic function of an unknown random vector in terms of a reference characteristic function}
\label{chGGC} 
The GGC derived as in Equation \eqref{GGC} can be used to give the characteristic function of an unknown \textit{pdf}, in terms of the characteristic function of a reference \textit{pdf}. For that taking Fourier transform $(\mathsf{F})$ of Equation \eqref{GGC}, we get:
\begin{align} 
\mathcal{F}_{\mathbf{x}}(\bs{\lambda}) &= \sum_{k=0}^{\infty} (-1)^k\frac{\bs{\alpha}(k,d)'}{k!}\mathsf{F}( \psi^{(k)}(\mathbf{x})) \\
\mathcal{F}_{\mathbf{x}}(\bs{\lambda}) & = \left[ \sum_{k=0}^{\infty}\bs{\alpha}(k,d)'  \frac{(i\bs{\lambda})^{\otimes k}}{k!}\right] \mathsf{F}(\psi(\mathbf{x})) \quad \mbox{($\because$ differentiation property of $\mathsf{F}$ )}\\
\label{GGCch}
& = \exp\left[ \sum_{k=1}^{\infty}\bs{\delta}(k,d)'  \frac{(i\bs{\lambda})^{\otimes k}}{k!}\right]\mathcal{F}_r(\mathbf{x}) \quad (\because \mbox{ Equation }\eqref{alphadeltaeq})
\end{align}
where, $\mathcal{F}_r$ is the characteristic function of the reference \textit{pdf}.
\section[Compact derivation for the GGC expansion]{Compact derivation for the Generalized Gram-Charlier expansion}
\label{compactGGC} 
The compact derivation of Equation \eqref{GGC} follows as under: 
\begin{align}
\mathcal{F}_{\mathbf{x}}(\bs{\lambda}) & = \exp\left[ \sum_{k=1}^{\infty}\mathbf{c}(k,d)'  \frac{(i\bs{\lambda})^{\otimes k}}{k!}\right] \quad \mbox{( $\because$ definition in Equation \eqref{ndCGF2} )}\\
& = \exp\left[ \sum_{k=1}^{\infty}\bs{\delta}(k,d)'  \frac{(i\bs{\lambda})^{\otimes k}}{k!}\right]\exp\left[ \sum_{k=1}^{\infty}\mathbf{c}_r(k,d)  \frac{(i\bs{\lambda})^{\otimes k}}{k!}\right] \\
& \quad (\because \bs{\delta}(k,d) = \mathbf{c}(k,d) - \mathbf{c}_r(k,d) ) \nonumber \\
& = \left[ \sum_{k=0}^{\infty}\bs{\alpha}(k,d)'  \frac{(i\bs{\lambda})^{\otimes k}}{k!}\right]\mathsf{F}(\psi(\mathbf{x}))& 
\end{align}
Taking inverse Fourier transform of the above equation brings 
\begin{align}
f(\mathbf{x}) &= \sum_{k=0}^{\infty}\bs{\alpha}(k,d)' \frac{(-1)^k}{k!}\bs{\delta}^{(k)}(\mathbf{x}) * \psi(\mathbf{x}) \\
\mbox{ or }f(\mathbf{x}) &= \sum_{k=0}^{\infty}\bs{\alpha}(k,d)' \frac{(-1)^k}{k!}\psi^{(k)}(\mathbf{x}) 
\end{align}
where, $*$ indicates convolution. Thus, the Equation \eqref{GGC} is obtained in a more compact way.
\section{Conclusion}
\label{conclusion}
The article has derived multivariate Generalized Gram-Charlier (GGC) expansion in Equation \eqref{GGC}; combined with Equation \eqref{ndalphadelta}; that expresses an unknown multivariate \textit{pdf} in terms of vector cumulants and vector derivatives of a reference \textit{pdf}. 
The multivariate Gram-Charlier A series is derived in Equation \eqref{ndGCAshort} and Equation \eqref{ndGCAhershort} representing an unknown multivariate \textit{pdf} in terms of  its vector cumulants and vector Hermite polynomials. There has been also derived compact formulas for obtaining multivariate  vector moments from vector cumulants  in Equation \eqref{ndmomcum} and vise-a-verse in Equation \eqref{ndcummom}; the integral form of  multivariate \textit{pdf} representation in Equation \eqref{ndpdfreal} and the integral form of multivariate vector Hermite polynomials in Equation \eqref{hermitenew}, as well, in Equation \eqref{hermitend01}.   The expressions are derived using only 
elementary calculus of several variables in vector notations through Kronecker product based derivative operator. Thus, they are  more transparent and more comprehensive compare to their corresponding multi-linear matrix representations or tensor representations. 
\appendix
\section{The Multivariate Representations of GCA Series and GGC Series }
\label{historyGGCreps}
As mention in Section \ref{introduction}, this section of the appendix describes existing representations of GCA series and GGC series. The goal is to place together various historical representations for the ease of comparison, on the level of difficulty or simplicity in representation, to the readers. 
Therefore, no attempt is made to explain their derivation or each terms in representation. For further details the actual references need be referred. 

\hspace{0.2 in} The GCA series representation using multi-element matrix notations for cumulants and moments by \citet{sauer1979convenient} is as under:
\begin{align}
\label{GCAsauer79}
f_{\mathbf{x}}(\mathbf{x}) &= \sum_{s_1=0}^{\infty} \sum_{s_2=0}^{\infty}\ldots \sum_{s_d=0}^{\infty} \left[ \mathbf{C}_{s_1s_2\cdots s_d} \cdot (-1)^{\sum_{i=1}^{d}s_i} \prod_{p=1}^{d} H_{s_p}(x_p)G(x_p) \right]  \\
& \quad \mbox{with, }  \mathbf{C}_{s_1s_2\cdots s_d} = \frac{E\left\{ \prod_{i=1}^{d} H_{s_i}(X_i) \right\} }{ (-1)^{\sum_{i=1}^{d}s_i} \prod_{j=1}^{d} s_j!}
\end{align}
where,  $\mathbf{C}_{s_1s_2\cdots s_d}$ is the constant depending upon cross-moments and $H_i(x)$ is the one-dimensional Hermite polynomial of $i^{th}$ order. 

\hspace{0.2 in}	The GCA Series representation using recursive formula for Hermite polynomials by \citet{ndHPGCberkowitz70} is as under:
\begin{align}
\label{GCAberkowitz70}
f_{\mathbf{x}}(\mathbf{x}) &= G(\mathbf{x})\sum_{m=0}^{\infty} A_m H_m(\mathbf{z})  \\
& \quad \mbox{with, } A_m = \prod_{i=1}^{d} (m_i!)^{-1} \int_{\mathbb{R^d}} J_m(z)f_\mathbf{x}(\mathbf{x}) d\mathbf{x} 
\end{align}
where, $G(\mathbf{x})$ denote the multivariate Gaussian; $\mathbf{z}=(\mathbf{x}-\bs{\mu})$ is the normalized variable; $\bs{\mu}$ is the mean vector; $\{ H_m(\mathbf{x})\}$ and $\{ J_m(\mathbf{x})\}$ are the complete bi-orthogonal system of Hermite polynomials. \\
Using recursive relations of $\{ J_m(\mathbf{x})\}$, $A_m$ is given recursively as under:
\begin{align}
A_m &= \frac{1}{N} \left[ \prod_{i=1}^{d} (m_i!)^{-1} \sum_{i=1}^{N}z_k^{(i)}J_{m-e_k}(y^{i}) - \sum_{i=1}^{N}r_{ki}m_k^{-1}A_{m-e_k-e_f} \right], k = 1,\ldots,d 
\end{align}
where, N is the number of available samples and $\mathbf{e}_k$ is a vector with a "1" as the $k^{th}$ component and "0" elsewhere. This defines the coefficients of expansions also recursively. 
  
\hspace{0.2 in}	The GGC series representation using tensor notations for cumulants and Hermite polynomials by \citet[Chapter 5]{Tensorbook87} is as under:
\begin{align}
\label{GGCMcCullagh87}
f_{\mathbf{x}}(x;\kappa) = f_0(\mathbf{x})\left[ 1+ \eta^ih_i(\mathbf{x}) + \eta^{ij} h_{ij}(\mathbf{x})/2! + \eta^{ijk} h_{ijk}(\mathbf{x})/3! + \ldots \right] 
\end{align} 
where, $h_i(\mathbf{x})= h_i(\mathbf{x};\bs{\lambda}) = f_i(\mathbf{x})/f_0(\mathbf{x}), \quad h_{ij}(\mathbf{x})= h_{ij}(\mathbf{x};\bs{\lambda}) = f_{ij}(\mathbf{x})/f_0(\mathbf{x}), \mbox{ }\ldots $ and 
$ f_i(\mathbf{x}) = \partial f_0(\mathbf{x})/\partial x^i, f_{ij}(\mathbf{x}) = \partial^2 f_0(\mathbf{x})/{\partial x^i \partial x^j}, \ldots$; so on. 
Also, given \\
$\kappa^i, \kappa^{i,j}, \kappa^{i,j,k}, \ldots $ are the cumulant tensors of random vector $\mathbf{x}$ and $\lambda^i, \lambda^{i,j}, \lambda^{i,j,k}, \ldots $ are the cumulant tensors of the reference \textit{pdf} $f_0(\mathbf{x})$; we get:
$$ \eta^i = \kappa^i - \lambda^i, \eta^{i,j} = \kappa^{i,j} - \lambda^{i,j}, \eta^{i,j,k} = \kappa^{i,j,k} - \lambda^{i,j,k}, \mbox{ }\ldots $$
The formal `moments'  $\eta^{i}, \eta^{ij}, \eta^{ijk}, \ldots $ are defined based on the formal `cumulants' (or the cumulant differences) $\eta^{i}, \eta^{i,j}, \eta^{i,j,k} $ and so on.  

\hspace{0.2 in}	Taking $f_0(\mathbf{x}) = G(\mathbf{x})$ i.e. multivariate Gaussian density as the reference \textit{pdf} and taking $\eta^i =0, \eta^{i,j} = 0$ in above Equation \eqref{GGCMcCullagh87};  the GCA series based on cumulant tensors is written as under:
\begin{align}
\label{GCAMcCullagh87}
f_{\mathbf{x}}(x;\kappa) &= G(\mathbf{x})\left[ 1+  \kappa^{i,j,k} h_{ijk}(\mathbf{x})/3! + \kappa^{i,j,k,l} h_{ijkl}(\mathbf{x})/4! + \kappa^{i,j,k,l,m} h_{ijklm}(\mathbf{x})/5! \right. \nonumber \\
& \quad \left.  + \left( \kappa^{i,j,k,l,m,n} + 10\kappa^{i,j,k}\kappa^{l,m,n} \right) h_{ijklmn}(\mathbf{x})/6!\ldots \right] 
\end{align} 
As could be observed, the GGC series and GCA series using tensor notations adds quite an ease to representation. But, with increase in number of terms, the difficulty in representation increases. 

\hspace{0.2 in}	The GCA series using vector moments and vector Hermite polynomials by \citet{dnHermite96} is as under:
\begin{align}
\label{GCAHermite96}
f_{\mathbf{x}}(\mathbf{x}) &= G(\mathbf{x}-\bs{\mu}; \mathbf{C}_{\mathbf{x}})\sum_{k=0}^{\infty} \frac{1}{k!} G_k^{T}(\mathbf{x}-\bs{\mu}; \mathbf{C}_{\mathbf{x}}) E\left\{ \mathbf{H}_k(\mathbf{x}-\bs{\mu}; \mathbf{C}_{\mathbf{x}}) \right\} 
\end{align}
where, $G_k(\mathbf{x} -\bs{\mu})$ is the $k^{th}$ order vector derivative of $G(\mathbf{x} -\bs{\mu})$ and 
$E\left\{ \mathbf{H}_k(\mathbf{x}-\bs{\mu}) \right\}$ is the expectation of $k^{th}$ order vector Hermite polynomial that is the function of vector moments. 
\begin{align}
E_f\left\{ \mathbf{H}_k(\mathbf{X}-\bs{\mu}; \mathbf{C}_\mathbf{x}^{-1})^{\otimes j}\right\} = k!\mathbf{S}_{d1_k} \sum_{j=0}^{[k/2]} \frac{\mathbf{m}_{k-2j}\otimes (- Vec\mbox{ }\mathbf{C}_{\mathbf{x}})^{\otimes j}}{(k-2j)!j!2^j}
\end{align}
where, $E_f\left\{ (\mathbf{X}-\bs{\mu})^{\otimes j}\right\} \equiv \mathbf{m}_j(\bs{\mu})$.

\hspace{0.2 in}	The GCA series representation using Bell polynomials is obtained  by 
\citet{withers2014dual}. Here, the Bell polynomials are represented through cumulant tensors.  With $r=(r_1,r_2, \ldots,r_d) \in \mathbb{I}^{d}$, $\mathbf{x}^r = x_1^{r_1}\ldots x_d^{r_d}$, $r!= r_1!r_2!\ldots r_d!$ and $|r| = r_1+ r_2 + \ldots + r_d$; the GCA series is as under:
\begin{align}
\label{GCABell}
f_{\mathbf{x}}(\mathbf{x})/G(\mathbf{x}) - 1 &= \sum_{|r|\geq 3}^{\infty} \mathbf{B}_rH_r(x,\mathbf{C}_{\mathbf{x}})/r! \\
\mbox{where, } \mathbf{B}_r &= \sum_{j=0}^{|r|} \mathbf{B}_{r,j} \\
\mathbf{B}_{r,j} &= \frac{1}{k_1(r-j)} \sum_{l=1}^{r-j} \binom{r}{j} \left( j+1-\frac{r+1}{l+1} \right) k_{l+1} \mathbf{B}_{r-l,j} \\
k_r &= \kappa( X_{r1}, X_{r2},\ldots, X_{rd})
\end{align}
where, $k_r$ is the $r^{th}$ order cumulant and $\kappa(\mathbf{x})$ is the moment.

\hspace{0.2 in}	Finally, the existing 
representations described in this Section \ref{historyGGCreps} can be compared with those derived in the current article. 
, can be compared with the GCA series derived  in Equation \eqref{ndGCAshort} and the GGC series derived in Equation \eqref{GGC} combined with Equation \eqref{alphadeltaeq}. 
More specifically, the GCA series in Equation \eqref{GCAsauer79} using multi-element array representation and the GGC series in Equation \eqref{GGCMcCullagh87} using tensor representation can be compared with the GCA series derived  in Equation \eqref{ndGCAshort} and the GGC series derived in Equation \eqref{GGC}, combined with Equation \eqref{alphadeltaeq}, in vector notations. The comparison  demonstrates the ease obtained in representation.  The same advantages are also obtained for other intermediate results in the article. 
\section{Some important properties of the K-derivative operator}
\label{KderProp}
Some important properties of the K-derivative operator are listed below. 
\begin{property}[Scaling Property] 
Let $\bs{\lambda} = (\lambda_1, \lambda_2, \ldots, \lambda_d)'$, $\bs{\lambda} \in \mathbb{R}^d$, $\mathbf{f}(\bs{\lambda}) \in \mathbb{R}^{m}$ and $\mathbf{f}_1(\bs{\lambda}) = \mathbf{A}\mathbf{f}(\bs{\lambda})$, where $\mathbf{A}$ is an $n \times m$ matrix. Then 
\begin{align}
\mathbf{D}_{\bs{\lambda}}^\otimes (\mathbf{f}_1) = \left( \mathbf{A} \otimes \mathbf{I}_d \right) \mathbf{D}_{\bs{\lambda}}^\otimes (\mathbf{f})
\end{align}
where, $\mathbf{I}_d$ is a d-dimensional unit matrix.
\end{property} 
\begin{property}[Chain Rule] 
Let $\bs{\lambda} = (\lambda_1, \lambda_2, \ldots, \lambda_d)'$, $\bs{\lambda} \in \mathbb{R}^d$, $\mathbf{f}(\bs{\lambda}) \in \mathbb{R}^{m_1}$  and $\mathbf{g}(\bs{\lambda}) \in \mathbb{R}^{m_2}$. Then 
\begin{align}
\mathbf{D}_{\bs{\lambda}}^\otimes (\mathbf{f} \otimes \mathbf{g}) = \mathbf{K}^{-1}_{3 \leftrightarrow 2}(m_1, m_2, d)((\mathbf{D}_{\bs{\lambda}}^\otimes \mathbf{f})\otimes \mathbf{g}) + \mathbf{f} \otimes (\mathbf{D}_{\bs{\lambda}}^\otimes \mathbf{g})
\end{align}
where, 
$\mathbf{K}_{3 \leftrightarrow 2}(m_1, m_2, d)$ is a commutation matrix of size $m_1m_2d \times m_1m_2d$ that changes the order of the the Kronecker product components. For example, 
\begin{align*}
\mathbf{K}_{3 \leftrightarrow 2}(m_1, m_2, d)(\mathbf{a_1}\otimes \mathbf{a_2} \otimes \mathbf{a_3} ) =  \mathbf{a_1} \otimes \mathbf{a_3} \otimes \mathbf{a_2}
\end{align*}
\end{property}
The above K-derivative properties can be used to derive further the following properties:
\begin{enumerate}
\item $ \mathbf{D}_{\bs{\lambda}}^\otimes \bs{\lambda} = Vec \mbox{ }\mathbf{I}_d $
\item $ \mathbf{D}_{\bs{\lambda}}^\otimes \bs{\lambda}^{\otimes k} = k \left( \bs{\lambda}^{\otimes k-1} \otimes Vec \mbox{ } \mathbf{I}_d \right)$
\item Let $g(\bs{\lambda}) = \mathbf{a}'f(\bs{\lambda})$ \\
$\mathbf{D}_{\bs{\lambda}}^\otimes g(\bs{\lambda}) = \left( \mathbf{a}\otimes \mathbf{I}_d \right)' \mathbf{D}_{\bs{\lambda}}^\otimes f(\bs{\lambda}) $  
\item $\mathbf{D}_{\bs{\lambda}}^\otimes \mathbf{a}'\bs{\lambda}^{\otimes k} = k \left( \mathbf{a}\otimes \mathbf{I}_d \right)' \left( \bs{\lambda}^{\otimes k-1} \otimes Vec \mbox{ } \mathbf{I}_d \right) = k \mbox{ }\mathbf{a}'\left( \bs{\lambda}^{\otimes k-1} \otimes \mathbf{I}_d \right)$\\
\item 
$\mathbf{D}_{\bs{\lambda}}^\otimes \mathbf{a}'^{\otimes k}\bs{\lambda}^{\otimes k} = k \left( \mathbf{a}^{\otimes k}\otimes \mathbf{I}_d \right)' \left( \bs{\lambda}^{\otimes k-1} \otimes Vec \mbox{ } \mathbf{I}_d \right) = k \left(\mathbf{a}' \bs{\lambda}\right)^{\otimes k-1}  \otimes \mathbf{a}$\\
\item $ \mathbf{D}_{\bs{\lambda}}^{\otimes r}\bs{\lambda}^{\otimes k} = k(k-1)\cdots (k-r+1) \left( \bs{\lambda}^{\otimes k-r} \otimes (Vec \mbox{ } \mathbf{I}_d)^{\otimes r} \right)$ \\
$\left( \mbox{ }\because \mbox{Repeated application of the Chain Rule } \right)$
\item 
$\mathbf{D}_{\bs{\lambda}}^{\otimes r} \mathbf{a}'^{\otimes k}\bs{\lambda}^{\otimes k} = k(k-1)\cdots (k-r+1) \left( \mathbf{a}^{\otimes k}\otimes \mathbf{I}_d^{\otimes r} \right)' \left( \bs{\lambda}^{\otimes k-r} \otimes (Vec \mbox{ } \mathbf{I}_d)^{\otimes r} \right) \\
\mbox{ }=  k(k-1)\cdots (k-r+1) \left(\mathbf{a}' \bs{\lambda}\right)^{\otimes k-r}  \otimes \mathbf{a}^{\otimes r}$\\
$\left( \mbox{ }\because \mbox{Repeated application of the Chain Rule and the property: } \right) $ 
\end{enumerate}
\section{The Taylor series expansion of some required functions near zero}
\label{Tayexpand}
The Taylor series expansion of the required functions near $\bs{\lambda} = \mathbf{0}$, based on the Equation  \eqref{ndTaylor}, are given as under: 
\begin{align}
\label{Tayeax}
\mbox{e}^{\mathbf{a}'\mathbf{x}} &= \sum_{k=0}^{\infty}\frac{ (\mathbf{a}'\mathbf{x})^{\otimes k}}{k!} = 1 + {\mathbf{a}'\mathbf{x}} + \frac{(\mathbf{a}'\mathbf{x})^{\otimes 2}}{2!} +  \frac{(\mathbf{a}'\mathbf{x})^{\otimes 3}}{3!} + \ldots \\
\label{Taysinx}
\sin({\mathbf{a}'\mathbf{x}}) &= \sum_{k=0}^{\infty}\frac{ (-1)^k }{(2k+1)!}(\mathbf{a}'\mathbf{x})^{\otimes 2k+1} = {\mathbf{a}'\mathbf{x}} - \frac{(\mathbf{a}'\mathbf{x})^{\otimes 3}}{3!} +  \frac{(\mathbf{a}'\mathbf{x})^{\otimes 5}}{5!} - \ldots \\
\label{Taycosx}
\cos(\mathbf{a}'\mathbf{x}) &= \sum_{k=0}^{\infty}\frac{ (-1)^k }{(2k)!}(\mathbf{a}'\mathbf{x})^{\otimes 2k} = 1 - \frac{(\mathbf{a}'\mathbf{x})^{\otimes 2}}{2!} +  \frac{(\mathbf{a}'\mathbf{x})^{\otimes 4}}{4!} - \ldots \\
\label{Tayeaxcosby}
\mbox{e}^{\mathbf{a}'\mathbf{x}}\cos(\mathbf{b}'\mathbf{y}) &= 1 + \mathbf{a}'\mathbf{x} + \frac{(\mathbf{a}'\mathbf{x})^{\otimes 2}}{2!} - \frac{(\mathbf{b}'\mathbf{y})^{\otimes 2}}{2!} +  \frac{(\mathbf{a}'\mathbf{x})^{\otimes 3}}{3!}  - \frac{(\mathbf{a}'\mathbf{x})\otimes (\mathbf{b}'\mathbf{y})^{\otimes 2}}{2!} \nonumber \\
& \mbox{ } + \frac{ (\mathbf{a}'\mathbf{x})^{\otimes 4} }{4!} + \frac{ (\mathbf{b}'\mathbf{y})^{\otimes 4} }{4!} -  \frac{(\mathbf{a}'\mathbf{x})^{\otimes 2}\otimes  (\mathbf{b}'\mathbf{y})^{\otimes 2} }{2!2!} + \ldots \\
\label{eaxsinby}
\mbox{e}^{\mathbf{a}'\mathbf{x}}\sin(\mathbf{b}'\mathbf{y}) &= \mathbf{b}'\mathbf{y}  + (\mathbf{a}'\mathbf{x})\otimes (\mathbf{b}'\mathbf{y})  - \frac{ (\mathbf{a}'\mathbf{x})^{\otimes 2}\otimes  (\mathbf{b}'\mathbf{y}) }{2!} 
- \frac{ ( \mathbf{b}'\mathbf{y} )^{\otimes 3} }{3!} \nonumber \\
& - \frac{ (\mathbf{a}'\mathbf{x}) \otimes  (\mathbf{b}'\mathbf{y})^{\otimes 3} }{3!}  +   
\frac{ ( \mathbf{b}'\mathbf{y} )^{\otimes 5} }{5!} -  \frac{ (\mathbf{a}'\mathbf{x})^{\otimes 2}\otimes  (\mathbf{b}'\mathbf{y})^{\otimes 3} }{2!3!}  + \ldots 
\end{align}
\section{Some Proofs}
\subsection{K-derivative of $\bs{\delta}(\mathbf{x})$}
\label{prfderdel}
Based on the Differentiation Property of Fourier Transform, the following is obtained:  
\begin{align}
 \mathsf{F}^{-1} \left( \mathsf{F}\left( \mathbf{D}^{\otimes k} f(\mathbf{x} \right) \right) 
& = \mathsf{F}^{-1} \left( (i\bs{\lambda})^{\otimes k} \mathsf{F}( f(\mathbf{x})) \right) \nonumber \\
  \Rightarrow \mathbf{D}^{\otimes k} \delta(\mathbf{x}) &= \mathsf{F}^{-1} ( i\bs{\lambda} )^{\otimes k}   \nonumber
\end{align} 
\subsection{Relationship between $\mathbf{m}(k,d)$ and $\mathbf{c}(k,d)$ for $k=2$ and $k=3$}
\label{Apmomcum}
For $k=2$; using Equation \eqref{momkd}, Equation \eqref{momcumeq} and the differentiation rules for Kronecker product in Appendix \ref{KderProp}; the following is derived:
\begin{align}
& \mathbf{m}(2,d) = \mathbf{D}_{\bs{\lambda}}^{\otimes 2} \mathbf{M}(\bs{\lambda})\vline_{\bs{\lambda}=\mathbf{0}} 	\\
& \Rightarrow 
\sum_{p=2}^{\infty} \mathbf{m}(p,d)'\frac{ \left( \bs{\lambda}^{\otimes (p-2)}\otimes \mathbf{I}_d \otimes \mathbf{I}_d \right)}{(p-2)!}\vline_{\bs{\lambda}= \mathbf{0}} = 
\mathbf{K}^{-1}_{\mathfrak{p}3\leftrightarrow 2}(d_1,d_2,d) \left\{ \left( \sum_{p=2}^{\infty}\mathbf{c}(p,d)' \right. \right. \nonumber \\
& \quad \left. \left. \frac{ \left( \bs{\lambda}^{\otimes (p-2)} \otimes \mathbf{I}_d \otimes \mathbf{I}_d \right) }{(p-2)!} \right)  
 \otimes \exp \left( \sum_{q=1}^{\infty}\mathbf{c}(q,d)'\frac{\bs{\lambda}^{\otimes q}}{q!}\right)\right\} \vline_{\bs{\lambda}= \mathbf{0}}	 \nonumber \\
& \quad +   \left\{ \left( \sum_{p=1}^{\infty}\mathbf{c}(p,d)'\frac{ \left( \bs{\lambda}^{\otimes (p-1)}\otimes \mathbf{I}_d \right) }{(p-1)!} \right)^{\otimes 2}  \otimes \exp \left( \sum_{q=1}^{\infty}\mathbf{c}(q,d)'\frac{\bs{\lambda}^{\otimes q}}{q!}\right)\right\}\vline_{\bs{\lambda}= \mathbf{0}}  \nonumber \\
\Rightarrow & \mathbf{m}(2,d) = \mathbf{c}(2,d) + \mathbf{c}(1,d)^{\otimes 2} 
\end{align}
where, $\mathfrak{p}3 \leftrightarrow 2 \in \mathfrak{P}$ is the required permutation and $\mathbf{K}_{\mathfrak{p}3\leftrightarrow 2}(d_1,d_2,d)$ is the corresponding commutation matrix.  \\
For $k=3$; using Equation \eqref{momkd},  Equation \eqref{momcumeq} and the differentiation rules for Kronecker product in Appendix \ref{KderProp}; the following is derived:
\begin{align}
& \mathbf{m}(3,d) = \mathbf{D}_{\bs{\lambda}}^{\otimes 3} \mathbf{M}(\bs{\lambda})\vline_{\bs{\lambda}=\mathbf{0}} 	\\
 \Rightarrow & \sum_{p=3}^{\infty} \mathbf{m}(p,d)'\frac{ \left( \bs{\lambda}^{\otimes (p-3)} \otimes \mathbf{I}_{d^3} \right) }{(p-3)!}\vline_{\bs{\lambda}= \mathbf{0}} = \mathbf{K}^{-1}_{\mathfrak{p}3\leftrightarrow 2}(d_1,d_2,d)\mathbf{K}^{-1}_{\mathfrak{p}3\leftrightarrow 2}(d_1,d_2,d) \nonumber \\ & \quad 
\left\{ \left( \sum_{p=3}^{\infty}\mathbf{c}(p,d)' 
\frac{ \left( \bs{\lambda}^{\otimes (p-3)}\otimes \mathbf{I}_{d^3} \right) }{(p-3)!} \right) \otimes \exp \left( \sum_{q=1}^{\infty}\mathbf{c}(q,d)'\frac{\bs{\lambda}^{\otimes q}}{q!}\right)\right\}\vline_{\bs{\lambda}= \mathbf{0}}	
\nonumber \\ & \quad 
+ \mathbf{K}^{-1}_{\mathfrak{p}3\leftrightarrow 2}(d_1,d_2,d) \left\{ \left( \sum_{p=2}^{\infty}\mathbf{c}(p,d)' 
\frac{ \left(  \bs{\lambda}^{\otimes (p-2)}\otimes \mathbf{I}_{d^2} \right) }{(p-2)!} \right) \otimes \left( \sum_{p=1}^{\infty}\mathbf{c}(p,d)' \right. \right. \nonumber \\
& \quad \left.  \left.
\frac{ \left( \bs{\lambda}^{\otimes (p-1)}\otimes \mathbf{I}_{d} \right) }{(p-1)!} \right) 
\otimes \exp \left( \sum_{q=1}^{\infty}\mathbf{c}(q,d)'\frac{\bs{\lambda}^{\otimes q}}{q!}\right)\right\}\vline_{\bs{\lambda}= \mathbf{0}} 
+ \mathbf{K}^{-1}_{\mathfrak{p}3\leftrightarrow 2}(d_1,d_2,d)  \nonumber \\
& \quad  
\left\{ 2 \left( \sum_{p=2}^{\infty}\mathbf{c}(p,d)'\frac{ \left(  \bs{\lambda}^{\otimes (p-2)}\otimes \mathbf{I}_{d^2} \right) }{(p-2)!} \right)^{\otimes 1}  
\otimes \left( \sum_{p=1}^{\infty}\mathbf{c}(p,d)' 
\frac{ \left( \bs{\lambda}^{\otimes (p-1)}\otimes \mathbf{I}_{d} \right) }{(p-1)!} \right)^{\otimes 1} 
\right. \nonumber \\ & \quad \left. 
\otimes  \exp \left( \sum_{q=1}^{\infty}\mathbf{c}(q,d)'\frac{\bs{\lambda}^{\otimes q}}{q!}\right)\right\}\vline_{\bs{\lambda}= \mathbf{0}}   
+ \left( \sum_{p=1}^{\infty}\mathbf{c}(p,d)' 
\frac{ \left( \bs{\lambda}^{\otimes (p-1)}\otimes \mathbf{I}_{d} \right) }{(p-1)!} \right)^{\otimes 2} 
\otimes \left( \sum_{p=1}^{\infty}\mathbf{c}(p,d)' \right. \nonumber \\& \quad \left. 
\frac{ \left( \bs{\lambda}^{\otimes (p-1)}\otimes \mathbf{I}_{d} \right) }{(p-1)!} \right)^{\otimes 1}  
\exp \left( \sum_{q=1}^{\infty}\mathbf{c}(q,d)'\frac{\bs{\lambda}^{\otimes q}}{q!}\right)\vline_{\bs{\lambda}= \mathbf{0}}  \nonumber \\
\Rightarrow & \mathbf{m}(3,d) = \mathbf{c}(3,d) + 3\mathbf{c}(2,d) \otimes \mathbf{c}(1,d) +\mathbf{c}(1,d)^{\otimes 3} 
\end{align}
\subsection{Multivariate Gaussian representation in terms of the cumulants (The proof related to the Section \ref{pdfcum}} 
\label{Gaussint}
Applying $\mathbf{D}_{\mathbf{x}}^{\otimes }$ on Equation \eqref{Gpdf}, 
\begin{align}
\label{dergx}
\mathbf{D}_{\mathbf{x}}^{\otimes }f(\mathbf{x}) = \frac{1}{(2\pi)^d}& \int_{\mathbb{R}^{d}} \left( \bs{\lambda} \otimes \mathbf{I}_d \right) \otimes 
\exp \left( - \frac{\mathbf{c}(2,d)'}{2}\bs{\lambda}^{\otimes 2} \right)   \sin\left( \left( \mathbf{x} - \mathbf{c}(1,d) \right)'\bs{\lambda} \right) d \bs{\lambda}  \\
\label{refeq1}
\mbox{Let, } 
U(\bs{\lambda}) &= \exp \left( - \frac{\mathbf{c}(2,d)'}{2}\bs{\lambda}^{\otimes 2} \right)  \\
\Rightarrow \mathbf{D}_{\bs{\lambda}}^{\otimes } U(\bs{\lambda}) &= - \mathbf{c}(2,d)'\left( \bs{\lambda} \otimes \mathbf{I}_d \right) \exp \left( - \frac{\mathbf{c}(2,d)'}{2}\bs{\lambda}^{\otimes 2} \right) \\
\label{refeq2}
\mbox{Also, let } V(\bs{\lambda}) &= \left( \mathbf{x}- \mathbf{c}(1,d) \right) \cos\left( \left( \mathbf{x}- \mathbf{c}(1,d) \right)'\bs{\lambda} \right) \\
\Rightarrow \int V d \bs{\lambda} &= \sin \left( \left( \mathbf{x}- \mathbf{c}(1,d) \right)'\bs{\lambda} \right) 
\end{align}
Using above Equation \eqref{refeq1}, Equation \eqref{refeq2} in Equation \eqref{dergx} and taking $ inv \mbox{ }\mathbf{c}(2,d) = Vec \left( \mathbf{C}_{\mathbf{x}}^{-1} \right)$; we get:  
\begin{align}
\label{eqdxfx4gx}
\mathbf{D}_{\mathbf{x}}^{\otimes }f(\mathbf{x}) &= - \frac{inv \mbox{ }\mathbf{c}(2,d)' }{(2\pi)^d}\int_{\mathbb{R}^d} \left( \mathbf{D}_{\bs{\lambda}}^{\otimes } U(\bs{\lambda}) \right) \left( \int V d \bs{\lambda} \right) d\bs{\lambda} 
\end{align}
Applying integration by parts to above Equation \eqref{eqdxfx4gx}, we get:
\begin{align}
\mathbf{D}_{\mathbf{x}}^{\otimes }f(\mathbf{x}) & = - \frac{inv \mbox{ }\mathbf{c}(2,d)' }{(2\pi)^d} \left\{ \left. \exp \left( -\frac{1}{2}\mathbf{c}(2,d)'\bs{\lambda}^{\otimes 2} \right) \sin\left( \left( \mathbf{x}- \mathbf{c}(1,d) \right)'\bs{\lambda} \right)  \right|_{\mathbb{R}^d} \right. \nonumber \\
& \quad \left. -  \exp \left( - \frac{\mathbf{c}(2,d)'}{2}\bs{\lambda}^{\otimes 2} \right) \left( \mathbf{x} - \mathbf{c}(1,d)\right)  \cos\left( \left( \mathbf{x} - \mathbf{c}(1,d)\right)'\bs{\lambda} \right) d \bs{\lambda}
\right\} \\ 
& = inv \mbox{ }\mathbf{c}(2,d)' \left( \mathbf{x} - \mathbf{c}(1,d)\right) f(\mathbf{x}) 
\end{align}
The solution of above differential equation leads to the following:
\begin{align}
f(\mathbf{x}) &= c \exp\left( - \left( inv \mbox{ }\mathbf{c}(2,d)\right)'( \mathbf{x}- \mathbf{c}(1,d) ) \right) \\
\mbox{where, } c &= f(\mathbf{0}) = \frac{1}{(2\pi)^{d}}\int_{\mathbb{R}^d}\exp \left( - \frac{\mathbf{c}(2,d)'}{2}\bs{\lambda}^{\otimes 2} \right)d\bs{\lambda} \nonumber \\
 & = (2\pi)^{-d/2}\left(|\mathbf{C}_{\mathbf{x}}| \right)^{-1/2}  \\
\Rightarrow  f(\mathbf{x}) &=   (2\pi)^{-d/2}\left(|\mathbf{C}_{\mathbf{x}}| \right)^{-1/2} 
\exp\left( - \left( inv \mbox{ }\mathbf{c}(2,d)\right)'(\mathbf{x}- \mathbf{c}(1,d)) \right) \\
\Rightarrow \frac{1}{(2\pi)^d} & \int_{\mathbb{R}^{d}}  \exp \left( - \frac{\mathbf{c}(2,d)'}{2}\bs{\lambda}^{\otimes 2} \right)  \cos\left( \mathbf{x}'\bs{\lambda} - \mathbf{c}(1,d)'\bs{\lambda} \right) d \bs{\lambda}  = G(\mathbf{x}) 
\end{align}
%
\bibliographystyle{elsarticle-harv}
\vskip 0.2in
\bibliography{Allref}

\end{document}